%% file: log.tex
 \documentclass[12pt]{article}

\usepackage{graphicx}
\usepackage{latexsym,amsmath,amsfonts,amscd, amsthm}
\usepackage{changebar}
\usepackage{color}
\usepackage{bm}
\usepackage{tikz}
\usepackage{multirow}
\usepackage{epsfig,subfigure}
\usepackage{multirow}
\usetikzlibrary{arrows,backgrounds,snakes}

\newtheoremstyle{plainNoItalics}{}{}{\normalfont}{}{\bfseries}{.}{ }{}

\theoremstyle{plain}
\newtheorem{thm}{Theorem}[section]

\theoremstyle{plainNoItalics}

\newtheorem{defn}[thm]{Definition}
\newtheorem{rem}[thm]{Remark}
\newtheorem{prop}[thm]{Proposition}
\newtheorem{exa}[thm]{Example}

\newcommand{\f}{\frac}

\newcommand{\beq}{\begin{equation}}
\newcommand{\eeq}{\end{equation}}
\newcommand{\beqa}{\begin{eqnarray}}
\newcommand{\eeqa}{\end{eqnarray}}
\newcommand{\bit}{\begin{itemize}}
\newcommand{\eit}{\end{itemize}}
\newcommand{\bedef}{\begin{defn}}
\newcommand{\edefn}{\end{defn}}
\newcommand{\bpro}{\begin{prop}}
\newcommand{\epro}{\end{prop}}

\begin{document}

\baselineskip=1.6pc


\input{title}

\newpage

\input{intro}
\input{RK_DG_MPP}
\input{simulations}
\input{conclusion}

\bibliographystyle{siam}
\bibliography{refer}

\end{document}

%% file: title.tex
\begin{center}
{\bf
High order maximum principle preserving discontinuous Galerkin method 
for convection-diffusion equations
}
\end{center}
\vspace{.2in}
\centerline{
Tao Xiong \footnote{Department of
Mathematics, University of Houston, Houston, 77204. E-mail:
txiong@math.uh.edu}
Jing-Mei Qiu \footnote{Department of Mathematics, University of Houston,
Houston, 77204. E-mail: jingqiu@math.uh.edu.
The first and second authors are supported by Air Force Office of Scientific Computing YIP grant FA9550-12-0318, NSF grant DMS-1217008, University of Houston.}
Zhengfu Xu \footnote{Department of
Mathematical Science, Michigan Technological University, Houghton, 49931. E-mail: zhengfux@mtu.edu. Supported by NSF grant DMS-1316662.}
}

\bigskip
\centerline{\bf Abstract}

{In this paper, we propose to apply the parametrized maximum-principle-preserving (MPP) flux
limiter in {\it{[Xiong et. al., JCP, 2013]}} to the discontinuous Galerkin (DG) method for solving
the convection-diffusion equations. The feasibility of applying the MPP flux limiters to the DG solution of convection-diffusion problem is based on the fact that the cell averages for the DG solutions are updated in a conservative fashion (by using flux difference) even in the presence of diffusion terms. The main purpose of this paper is to address the difficulty of obtaining higher than second order accuracy while maintaining a discrete maximum principle for the DG method solving convection diffusion equations. We found that the proposed MPP flux limiter can be applied to arbitrarily high order DG method. Numerical evidence is presented to show that the proposed MPP flux limiter method does not adversely affect the desired high order accuracy, nor does it require restrictive time steps. Numerical experiments including
incompressible Navier-Stokes equations demonstrate the high order accuracy preserving, the MPP performance,
and the robustness of the proposed method.}
\vfill

\newpage


%% file: intro.tex
\section{Introduction}
\label{sec1}
\setcounter{equation}{0}
\setcounter{figure}{0}
\setcounter{table}{0}

In this paper, we propose a parametrized maximum-principle-preserving (MPP) flux limiter
for the high order discontinuous Galerkin (DG) finite element method in order to solve the nonlinear convection-diffusion equation 
\begin{equation}
u_t+f(u)_x=a(u)_{xx}, \qquad u(x,0)=u_0(x).
\label{ad1}
\end{equation}
The exact solution of (\ref{ad1}) satisfies the maximum principle, that is, if
\begin{equation}
u_M=\max_x u_0(x), \qquad u_m=\min_x u_0(x),
\end{equation}
we have 
\begin{equation}
\label{CMP}
u(x,t)\in [u_m, u_M], \quad \forall t>0.
\end{equation} 
When $u(x, t)$ describes the density of a particular species, $u_m=0$, the problem is generally addressed as positivity preserving.



High order shock-capturing numerical methods for convection-dominated problems include the high resolution finite volume (FV) and finite difference (FD) essentially non-oscillatory (ENO) and weighted ENO (WENO) methods for the convection part, which are capable of producing solutions with fidelity without spurious oscillations. In this framework, high order central difference is generally used to approximate the second order derivative terms. See the lecture notes \cite{shu1998essentially,shu1999high} and the review paper \cite{shu2009high} of Shu and reference therein for more discussion of numerical methods in this aspect. The discontinuous Galerkin (DG) method in the finite element framework is another type method; it was developed by Cockburn et. al. in a series of papers \cite{cockburn1989tvb2,cockburn1989tvb3,cockburn1990runge,cockburn1998runge} for hyperbolic conservation laws and systems. The DG method has been well-known for its flexibility, h-p adaptivity, compactness and high parallel efficiency~\cite{cockburn2000development}. Later the DG method was generalized to the convection-dominated diffusion equations. Different types of DG approaches for solving the convection diffusion equations include the local DG (LDG) method \cite{cockburn1998runge}, the DG formulation of Cheng and Shu \cite{cheng2008discontinuous}. 
When a convection-dominated diffusion problem is solved within any of the two previously mentioned frameworks, numerical solutions may exhibit overshoots or undershoots, i.e, a discrete version of the maximum principle
\begin{equation}
\label{DMP}
u_m \le u^n_j \le u_M, \forall n, j,
\end{equation}
is no longer satisfied.

 
Two kinds of high order discrete maximum-principle-preserving limiters are newly developed for convection dominated problems. One is the polynomial rescaling MPP limiter proposed by Zhang and Shu in \cite{zhang2010maximum,zhang2011maximum} for hyperbolic conservation laws. It has been extended to the convection-diffusion equations based on a twice-integrated FV formulation of (\ref{ad1}) within the FV high order WENO framework \cite{zhang2012maximumcd}. The same technique under the DG framework for hyperbolic conservation laws has been applied to the convection-diffusion equations on a triangular mesh in \cite{yzhang2012maximum}, however the approach does not work for the Runge-Kutta DG (RKDG) method with order higher than $2$. The high order parametrized MPP flux limiter is developed in \cite{mpp_xu, mpp_xuMD} for hyperbolic conservation laws, which was later improved by Xiong et. al. \cite{mpp_xqx} by applying the limiter only at the final stage of a high order RK method. The MPP flux limiter has been generalized to convection-diffusion equations under the FD WENO framework in \cite{jiang2013parametrized} and the FV WENO framework in \cite{yang2014parametrized}. Early discussion of the discrete maximum principle for the convection diffusion equations includes the linear finite element solutions for parabolic equations \cite{fujii1973some} with recent developments in \cite{farago2005discrete,farago2006discrete,farago2012discrete,vejchodsky2010discrete} and the Petrov-Galerkin finite element method for convection-dominated problems \cite{mizukami1985petrov}.
However, they are under a different framework. 

In this paper, we propose to apply the parametrized MPP flux limiter in \cite{mpp_xqx} to the RKDG method, for solving convection-diffusion equations. For the convection part, the parametrized MPP flux limiter is proposed to preserve the MPP property {for the cell averages and at the final RK stage only}. This is different from the polynomial rescaling limiter proposed by Zhang and Shu that preserves the MPP property for the entire polynomial (or at Gaussian quadrature points) per element and at each of the RK stages. For the diffusion part, the parametrized MPP flux limiter is proposed for the DG formulations in \cite{cheng2008discontinuous, cockburn1998runge} with general piecewise $P^k$ ($k\ge0$) polynomial solution spaces. By the design, the parametrized MPP flux limiter preserves the MPP property of cell averages, thus avoids the main difficulty in the approach of polynomial rescaling limiter. Specifically, in \cite{yzhang2012maximum}, great effort is made to rewrite the updated solution as a convex combination of the solution point values in the current time step via approximating the second derivative term by point values; because of such complications, the limiter in \cite{yzhang2012maximum} is proposed for the DG schemes with $P^0$ and $P^1$ solution spaces only. Our proposed approach can be viewed as a low-cost easy-to-implement post-processing procedure that modifies the high order numerical fluxes towards a first order one {\em only at the final RK stage for the evolution of cell averages (not for higher moments)}, in order to preserve the solution cell averages' MPP property. We remark that the proposed DG solutions (piecewise polynomials) with the parametrized flux limiters might be out of the bound of $[u_m, u_M]$ with the cell averages well bounded by $[u_m, u_M]$. One can apply the polynomial rescaling limiters as in \cite{zhang2010maximum} {\em in the final time step only} to ensure the numerical solution (piecewise polynomials) to be within bounds. 

The proposed MPP flux limiter in the DG framework is mass conservative due to the flux difference form. It is very efficient due to the fact that the limiter is only applied at the final RK stage and for the cell averages per time step. The parametrized MPP flux limiter for the DG method can be proved to maintain up to 3rd order accuracy under the time step constraint of the original DG method, by following a similar analysis as in \cite{yang2014parametrized, jiang2013parametrized} for the finite volume and finite difference WENO method. The proof for higher than third order case is very technical and algebraically complicated, thus we rely on extensive numerical tests to showcase that up to fourth order accuracy can be preserved.  Extensive numerical tests, including the incompressible Navier-Stokes system, are presented to demonstrate the robust performance of the proposed approach in preserving the high order accuracy as well as the MPP property. 

The outline of the paper is as follows. In Section \ref{sec2}, the DG formulation of Cheng and Shu \cite{cheng2008discontinuous} for one and two dimensions are described. The application of the parametrized MPP flux limiters on the cell averages of the DG solution is presented. Extension to the LDG method will be discussed. Numerical results are provided in Section \ref{sec3}. Finally conclusions are made in Section \ref{sec4}.

%% file: RK_DG_MPP.tex
\section{The MPP flux limiter for the RKDG method}
\label{sec2}
\setcounter{equation}{0}
\setcounter{figure}{0}
\setcounter{table}{0}

In this section, we will first briefly describe the DG method developed by Cheng and Shu \cite{cheng2008discontinuous} for directly solving convection-diffusion equations.
Then we will apply the parametrized MPP flux limiters developed in \cite{mpp_xqx}
to the RKDG method. Both one and two dimensional cases will be presented.

\subsection{One dimensional case}
\label{sec2.1}

Without loss of generality, we assume periodic boundary condition or zero boundary condition with compact support for 1D cases. The spatial domain $[a,b]$ is discretized by $N$ cells,
\beq
a=x_{\f12}<x_{\f32} <\cdots <x_{N-\f12}<x_{N+\f12}=b,
\eeq
with the cell, cell center to be 
\beq
I_j=(x_{j-\f12}, x_{j+\f12}),\qquad x_j=\f12(x_{j-\f12}+x_{j+\f12}), \quad \forall j = 1,\cdots N.
\eeq 
The mesh size $h_j=x_{j+\f12}-x_{j-\f12}$ and let $h=\max_j h_j$. In the following, for simplicity, we
assume uniform mesh sizes $h_j=h, \forall j$. 

The DG method in \cite{cheng2008discontinuous} is defined as follows:
find $u_h\in V^k_h$, such that
\beqa
\int_{I_j}(u_h)_tv_hdx&-&\int_{I_j}f(u_h)(v_h)_xdx-\int_{I_j}a(u_h)(v_h)_{xx}dx \nonumber \\
&+&(\hat{f}(u^-_h,u^+_h)v^-_h)_{j+\f{1}{2}}-(\widehat{f}(u^-_h,u^+_h)v^+_h)_{j-\f{1}{2}}
-(\widetilde{a(u_h)_x}v^-_h)_{j+\f{1}{2}}  \nonumber \\
&+&(\widetilde{a(u_h)_x}v^+_h)_{j-\f{1}{2}}+(\widehat{a(u_h)}(v_h)^-_x)_{j+\f{1}{2}}
-(\widehat{a(u_h)}(v_h)^+_x)_{j-\f{1}{2}}  \nonumber \\
&=& 0, \label{dg1}
\eeqa
for any test function $v_h \in V^k_h$ and $j=1, \dots, N$, where $V^k_h=\{v: v|_{I_j}\in P^k(I_j),\forall j\}$ and $P^k(I_j)$ is a piecewise polynomial space with degree up to $k$ on the cell $I_j$. $\widehat{f}(u^-_h,u^+_h)$ is a monotone flux for the convection term, $\widetilde{a(u_h)_x}$ and $\widehat{a(u_h)}$
are numerical fluxes chosen to be
\beq
\widetilde{a(u_h)_x}=\f{[a(u_h)]}{[u_h]}((u_h)^-_x+\f{\alpha}{h}[a(u_h)]),\qquad \widehat{a(u_h)}=a(u^+_h),
\label{flux1}
\eeq
here $\alpha$ is a positive constant chosen for stability. $\xi_{j+\frac12}^-$ and $\xi_{j+\frac12}^+$ are the left and right limit from cell $I_j$ and $I_{j+1}$ respectively. $[\xi]_{j+\frac12}=\xi_{j+\frac12}^+ - \xi_{j+\frac12}^-$ is the jump of $\xi$ at the cell interface $x_{j+\frac12}$. A third order strong stability preserving (SSP) RK time discretization \cite{shu1988efficient} for the semi-discrete scheme (\ref{dg1}) is
given as
\begin{eqnarray}
\nonumber u^{(1)} &=& u^n+\Delta t L(u^n),\\
          u^{(2)} &=& u^n+\Delta t(\frac{1}{4}L(u^{n})+\frac{1}{4}L(u^{(1)})),\label{eq: rk3}\\
\nonumber u^{n+1} &=& u^n+\Delta t (\frac{1}{6} L(u^{n})+\frac{2}{3} L(u^{(2)})+\frac{1}{6} L(u^{(1)})).
\end{eqnarray}

The parametrized MPP flux limiters are applied to keep {\em only} the cell averages of $u_h$ within the range of 
$[u_m, u_M]$. If we take $v_h=1$ in (\ref{dg1}) and divided by $h$ on both sides, we have
\beq
\f{d}{dt}\bar{u}_h+\f{1}{h}\left(\hat{H}_{j+\f{1}{2}}-\hat{H}_{j-\f{1}{2}}\right)=0,
\label{cellaverage}
\eeq
where the flux $\hat{H}=\hat{f}(u^-_h,u^+_h)-\widetilde{a(u_h)_x}$. With the third order SPP RK method \eqref{eq: rk3}, 
the update of cell averages in equation \eqref{cellaverage} can be written as
\begin{eqnarray}
\label{eq: RKall}
(\bar{u}_h)^{n+1}_j=(\bar{u}_h)^{n}_j-\lambda (\hat H^{rk}_{j+\f12}-\hat H^{rk}_{j-\f12}),
\end{eqnarray}
where
\beq
\hat H^{rk}_{j+\f12} \doteq \frac{1}{6}\hat H^n_{j+\f12} + \frac23 \hat H^{(2)}_{j+\f12}+\frac{1}{6}\hat H^{(1)}_{j+\f12},
\label{rkfinal}
\eeq
with $\hat H^*$ being the numerical flux obtained by $u_h^*$ at each RK stage for $*=n, (1), (2)$, respectively.

The MPP flux limiter is proposed to replace the numerical flux $\hat H^{rk}_{j+\f12}$ by a modified one
\begin{eqnarray}
\label {eq: rkfl}
\tilde H^{rk}_{j+\f12} =\theta_{j+\f12} (\hat H^{rk}_{j+\f12}-\hat h_{j+\f12})+\hat h_{j+\f12},
\end{eqnarray}
where $\hat h_{j+\f12}$ is a first order monotone flux with which the scheme is maximum principle preserving,
e.g., the global Lax-Friedrichs flux \cite{cockburn2001runge}. 
The parameter $\theta_{j+\f12}$ is defined to ensure $(\bar{u}_h)^{n+1}_j \in [u_m, u_M]$, for which 
sufficient inequalities to have are
\begin{eqnarray}
\label{eq: umax}
\lambda \theta_{j-\f12} (\hat H^{rk}_{j-\f12}-\hat h_{j-\f12}) - \lambda \theta_{j+\f12} (\hat H^{rk}_{j+\f12}-\hat h_{j+\f12})-\Gamma^M_j &\le& 0, \\
\label{eq: umin}
\lambda \theta_{j-\f12} (\hat H^{rk}_{j-\f12}-\hat h_{j-\f12}) - \lambda \theta_{j+\f12} (\hat H^{rk}_{j+\f12}-\hat h_{j+\f12})-\Gamma^m_j &\ge& 0,
\end{eqnarray}
with $\lambda=\Delta t/h$ and
\[
\Gamma^M_j=u_{M}-(\bar u_h)^n_j+\lambda (\hat h_{j+\f12}-\hat h_{j-\f12}) \ge 0,
\quad
\Gamma^m_j=u_{m}-(\bar u_h)^n_j+\lambda (\hat h_{j+\f12}-\hat h_{j-\f12}) \le 0.
\]
Let $F_{j\pm\f12}\doteq\hat H^{rk}_{j\pm\f12}-\hat h_{j\pm\f12}$, 
the parameter $\theta_{j+\f12}$ can be obtained as follows, for details see \cite{mpp_xqx}:
\begin{enumerate}
\item
Assume
$
\theta_{j-\f12} \in [0, \Lambda^M_{-\f12, I_j}], \quad
\theta_{j+\f12} \in [0, \Lambda^M_{+\f12, I_j}],
$
where $\Lambda^M_{-\f12, I_j}$ and $\Lambda^M_{+\f12, I_j}$
are designed to preserve the upper bound by equation \eqref{eq: umax},
\begin{enumerate}
\item If $F_{j-\f12}\le 0$ and $F_{j+\f12}\ge 0$,
$(\Lambda^M_{-\f12, I_j}, \Lambda^M_{+\f12, I_j})=(1, 1).
$
\item  If $F_{j-\f12}\le 0$ and $F_{j+\f12}< 0$,
$
(\Lambda^M_{-\f12, {I_j}}, \Lambda^M_{+\f12, {I_j}})=(1, \min(1, \frac{\Gamma^M_j}{-\lambda F_{j+\f12}})).
$
\item  If $F_{j-\f12}> 0$ and $F_{j+\f12}\ge 0$,
$
(\Lambda^M_{-\f12, {I_j}}, \Lambda^M_{+\f12, {I_j}})=(\min(1, \frac{\Gamma^M_j}{\lambda F_{j-\f12}}), 1).
$
\item  If $F_{j-\f12}> 0$ and $F_{j+\f12}< 0$,
\bit
\item
If equation  \eqref{eq: umax} is satisfied with $(\theta_{j-\f12}, \theta_{j+\f12})=(1, 1)$, then
$
(\Lambda^M_{-\f12, {I_j}}, \Lambda^M_{+\f12, {I_j}})=(1, 1).
$
\item
If equation \eqref{eq: umax} is not satisfied with $(\theta_{j-\f12}, \theta_{j+\f12})=(1, 1)$, then
$
(\Lambda^M_{-\f12, {I_j}}, \Lambda^M_{+\f12, {I_j}})=(\frac{\Gamma^M_j}{\lambda F_{j-\f12}- \lambda F_{j+\f12}},\frac{\Gamma^M_j}{\lambda F_{j-\f12}- \lambda F_{j+\f12}} ).
$
\eit
\end{enumerate}
\item Similarly
assume
\[
\theta_{j-\f12} \in [0, \Lambda^m_{-\f12, I_j}], \quad
\theta_{j+\f12} \in [0, \Lambda^m_{+\f12, I_j}],
\]
where $\Lambda^m_{-\f12, I_j}$ and $\Lambda^m_{+\f12, I_j}$
are designed to preserve the lower bound by equation \eqref{eq: umin},
\begin{enumerate}
\item If $F_{j-\f12}\ge 0$ and $F_{j+\f12}\le 0$,
$
(\Lambda^m_{-\f12, I_j}, \Lambda^m_{+\f12, I_j})=(1, 1).
$
\item  If $F_{j-\f12}\ge 0$ and $F_{j+\f12}> 0$,
$
(\Lambda^m_{-\f12, {I_j}}, \Lambda^m_{+\f12, {I_j}})=(1, \min(1, \frac{\Gamma^m_j}{-\lambda F_{j+\f12}})).
$
\item  If $F_{j-\f12}< 0$ and $F_{j+\f12}\le 0$,
$
(\Lambda^m_{-\f12, {I_j}}, \Lambda^m_{+\f12, {I_j}})=(\min(1, \frac{\Gamma^m_j}{\lambda F_{j-\f12}}), 1).
$
\item  If $F_{j-\f12}< 0$ and $F_{j+\f12}> 0$,
\bit
\item
If equation  \eqref{eq: umin} is satisfied with $(\theta_{j-\f12}, \theta_{j+\f12})=(1, 1)$, then
$
(\Lambda^m_{-\f12, {I_j}}, \Lambda^m_{+\f12, {I_j}})=(1, 1).
$
\item If equation  \eqref{eq: umin} is not satisfied with $(\theta_{j-\f12}, \theta_{j+\f12})=(1, 1)$, then
$
(\Lambda^m_{-\f12, {I_j}}, \Lambda^m_{+\f12, {I_j}})=(\frac{\Gamma^m_j}{\lambda F_{j-\f12}- \lambda F_{j+\f12}},\frac{\Gamma^m_j}{\lambda F_{j-\f12}- \lambda F_{j+\f12}} ).
$
\eit
\end{enumerate}
\end{enumerate}
The local parameter $\theta_{j+\f12}$ is determined to be
\begin{eqnarray}
\label{eq: limit}
\theta_{j+\f12}=\min(\Lambda^M_{+\f12, {I_j}}, \Lambda^M_{-\f12, {I_{j+1}}}, \Lambda^m_{+\f12, {I_j}}, \Lambda^m_{-\f12, {I_{j+1}}}),
\end{eqnarray}
by the consideration to ensure both the upper bound (\ref{eq: umax}) and lower bound (\ref{eq: umin}) of
the cell averages in both cell $I_j$ and $I_{j+1}$.

\begin{rem}
In the proposed approach, the convection and diffusion terms are treated together for the MPP property of the cell averages. The parametrized flux limiters are applied {\em only for cell averages (not for higher moments) and at the final RK stage as in \eqref{eq: RKall}}. Hence, the proposed flux limiting procedure has low computational cost and is easy to implement. The approach can also be generalized to other multi-stage RK and multi-step methods. 
We remark that the DG solutions (piecewise polynomials) with the parametrized flux limiters might be out of the bound of $[u_m, u_M]$ with the cell averages well bounded by $[u_m, u_M]$. One can apply the polynomial rescaling limiters as in \cite{zhang2010maximum} {\em in the final time step only} to ensure the numerical solution (piecewise polynomials) to be within bounds. 
\end{rem}
\begin{rem}
The proposed flux limiter is different from the polynomial rescaling techniques introduced in \cite{zhang2010maximum, yzhang2012maximum} in several aspects. First of all, in \cite{zhang2010maximum}, the entire polynomial (or at least the Gaussian-Lobatto quadrature points) over each cell at each of the RK stage are rescaled to satisfy the MPP property. As a result, the temporal accuracy for a multi-stage RK method may be affected, e.g. see discussions in \cite{zhang2010maximum}. Secondly, in \cite{yzhang2012maximum} the convection and diffusion terms are treated separately; for both terms, great effort has been made to rewrite the updated cell average as a convex combination of point values in the current time step. Such approach introduces extra CFL time step constraint on the DG method; moreover, it is difficult to generalize such approach for the diffusion term with higher than second order accuracy, see Remark 2.2 in \cite{yzhang2012maximum}. 
\end{rem}

\begin{rem} (accuracy)
For the DG method (\ref{dg1}) with the numerical fluxes (\ref{flux1}), it has been proved in \cite{cheng2008discontinuous} that the method is stable with sub-optimal error estimate ($k$-th order for $P^k$ polynomial space) for the $L^2$ norm. Numerically, both $L^1$ and $L^\infty$ norms are of the optimal $(k+1)$-th order. Regarding the preservation of high order accuracy of the original RK DG method with the proposed MPP flux limiter, similar conclusions can be established following the same line of proof in \cite{yang2014parametrized, jiang2013parametrized}, i.e. without any additional time step condition, (1) for the special case of linear advection problem, the high order accuracy of the original RK DG solutions is maintained with the proposed flux limiter; (2) for the general convection-dominated diffusion problem, up to third order accuracy will be maintained with the flux limiter. In fact, numerically, it can be shown that arbitrary high order accuracy is preserved, under the time step constraint from the linear stability analysis for the DG method \cite{cockburn2001runge,wang2013error}.  
\end{rem}

\begin{rem}
The parametrized MPP flux limiter can also be applied to the local DG (LDG) method
for the convection-diffusion equations \cite{cockburn1998local}. To obtain an LDG formulation
for (\ref{ad1}), first we rewrite it as
\begin{equation}
u_t+f(u)_x=(\gamma(u)q)_x, \quad q-\Gamma(u)_x=0,
\label{ldg0}
\end{equation}
where $\gamma(u)=\sqrt{a'(u)}$ and $\Gamma(u)=\int^u \gamma(s)ds$. The LDG method is defined to be:
find $u_h, q_h \in V^k_h$, such that for any test functions $v_h, w_h \in V^k_h$, we have
\begin{eqnarray}
\int_{I_j}(u_h)_t v_h dx-\int_{I_j}(f(u_h)-\gamma(u_h)q_h) (v_h)_x dx+(\hat f-\hat{\gamma}\hat q)_{j+\frac12}(v_h)^-_{j+\frac12}-(\hat f-\hat{\gamma}\hat q)_{j-\frac12}(v_h)^+_{j-\frac12}=0, \nonumber \\
\label{ldg1}\\
\int_{I_j}q_hw_hdx+\int_{I_j}\Gamma(u_h)(w_h)_xdx-\hat \Gamma_{j+\frac12}(w_h)^-_{j+\frac12}+\hat \Gamma_{j-\frac12}(w_h)^+_{j-\frac12}=0. \nonumber \\
\label{ldg2}
\end{eqnarray}  
$\hat f =\hat f(u^+_h,u^-_h)$ is the monotone flux for the convection part. For the diffusion part,
the numerical fluxes are
\begin{equation}
\hat \gamma=\frac{\Gamma(u^+_h)-\Gamma(u^-_h)}{u^+_h-u^-_h},\quad \hat q=q^-_h,\quad \hat \Gamma=\Gamma(u^+_h).
\end{equation}
If we take $v_h=1$ in (\ref{ldg1}), we have the same equation (\ref{cellaverage}) for the cell
average of $u_h$, the only difference is the flux given by $\hat H=\hat f-\hat \gamma \hat q$.
The rest of applying the flux limiter would be the same as described above. Similar arguments hold
for the two dimensional case in the following subsection.
\end{rem}

\begin{rem}
For convection-diffusion equations with source terms $u_t+f(u)_x=a(u)_{xx}+s(u,x)$, the technique
in \cite{pp_euler} can be used for the source term to ensure the MPP property.
\end{rem}

\subsection{Two dimensional case}
In this subsection, we consider the generalization of the parametrized flux limiter to the two dimensional convection-diffusion equation 
\begin{equation}
u_t + \nabla \cdot \mathbf{F}(u) = \nabla \cdot (\mathbf{A} \nabla u), \qquad \mathbf{F}(u)=(f(u),g(u)),
\label{ad2}
\end{equation}
on a bounded domain of $\mathbf{x}=(x,y)\in[a, b]\times[c, d]$, where $\mathbf{A}=\mathbf{A}(u,\mathbf{x})$ is a $2\times2$ symmetric semi-positive-definite matrix. Similar observation of (\ref{CMP}) also holds for the two dimensional case.

For simplicity, in the following, we assume periodic boundary conditions or zero boundary conditions with compact support in each direction. A spatial discretization with $N_x\times N_y$ rectangular meshes is defined as
\[
a=x_{\f12}<x_{\f32}<\cdots<x_{N_x-\f12}<x_{N_x+\f12}=b, \quad  c=y_{\f12}<y_{\f32}<\cdots<y_{N_y-\f12}<y_{N_y+\f12}=d,
\]
where the cell, cell centers and cell sizes are defined by
\[
K_{ij}=I_i\times J_j, \quad I_i=(x_{i+\f12},x_{i-\f12}), \quad J_j=(y_{j-\f12}, y_{j+\f12}), 
\]
\[
h^x_i=x_{i+\f12}-x_{i-\f12}, \quad x_i=\f12(x_{i+\f12}+x_{i-\f12}), \quad h^y_j=y_{j+\f12}-y_{j-\f12}, \quad y_j=\f12(y_{j-\f12}+y_{j+\f12}),  
\]
and $h^x=\max_{i}h^x_i$, $h^y=\max_jh^y_j$, $h=\max(h^x, h^y)$. For simplicity, in the following, we assume $h^x_i=h^x, \forall i$ and $h^y_j=h^y, \forall j$. 

The DG scheme in \cite{cheng2008discontinuous} for two dimensions with rectangular mesh is defined as: find $u_h\in V^k_h$, such that for any test function $v_h\in V^k_h$,
\beqa
\int_{K_{ij}}&(u_h)_t& v_hd\mathbf{x}-\int_{K_{ij}}\mathbf{F}(u_h) \cdot \nabla v_hd\mathbf{x} -\int_{K_{ij}}u_h\nabla\cdot{(\mathbf{A}\nabla v_h)}d\mathbf{x} \nonumber \\
&+&\int_{\partial K_{ij}}(\widetilde{u_h\mathbf{A}} \nabla v_h) \cdot \mathbf{n} ds 
+\int_{\partial K_{ij}}\left( \hat{\mathbf{F}}\cdot \mathbf{n}-\widehat{\mathbf{A}\nabla u_h} \cdot \mathbf{n} \right) v_h ds = 0. \label{dg2}
\eeqa
Here $V^k_h=\{v:v_{K_{ij}} \in P^k(K_{ij}),\forall i,j \}$ and $P^k(K_{ij})$ is the two dimensional polynomial space with degree up to $k$ on the cell $K_{ij}$, $\mathbf{n}$ is the outward unit normal vector
on the edges. $\hat{\mathbf{F}}=\hat{\mathbf{F}} (u^+_h,u^-_h)$
is a monotone numerical flux for the convection part \cite{cockburn1990runge}, e.g., the global Lax-Friedrichs flux. Other numerical fluxes are defined by \cite{yzhang2012maximum}
\beq
\label{flux2d}
\widehat{\mathbf{A}\nabla u_h} \cdot \mathbf{n}=\mathbf{A}(u^-_h)\nabla u^-_h \cdot \mathbf{n}
+\frac{\alpha \Lambda}{h}(u^{out}_h-u^{in}_h), \qquad \widetilde{u_h \mathbf{A}}=u^+_h\mathbf{A}(u^+_h),
\eeq
here $\Lambda$ is the maximum absolute eigenvalue of the symmetric matrix $\mathbf{A}$, $\alpha$ is a
parameter large enough to ensure the stability of the scheme, which will be specified later. $u^\pm_h$ are
the left and right limit values from the cells adjacent to the edges respectively. On the left boundary of $K_{ij}$, we have $u^{out}_h=u^-_h$ and $u^{in}_h=u^+_h$, while on the right boundary, $u^{out}_h=u^+_h$ and $u^{in}_h=u^-_h$. Similarly for $u^{out}_h$ and $u^{in}_h$ on the top and bottom boundaries of $K_{ij}$.

Taking $v_h=1$ in (\ref{dg2}), for the cell average, simply we have 
\beqa
\f{d}{dt}\bar{u}_h
&+&\f{1}{h^x}\left(\frac{1}{h^y}\int_{J_j}\hat{H}(x_{i+\f{1}{2}},y)dy
- \frac{1}{h^y}\int_{J_j}\hat{H}(x_{i-\f{1}{2}},y)dy\right) \nonumber \\
&+&\f{1}{h^y}\left(\frac{1}{h^x}\int_{I_i}\hat{G}(x,y_{j+\f{1}{2}})dx-\frac{1}{h^x}\int_{I_i}\hat{G}(x,y_{j-\f{1}{2}})dx\right)=0.
\label{cellaverage2}
\eeqa
where $\hat{H}$ and $\hat{G}$ are $\hat{\mathbf{F}}\cdot \mathbf{n}-\widehat{\mathbf{A}(\nabla u_h)} \cdot \mathbf{n}$ with $\mathbf{n}=(1,0)$ and $\mathbf{n}=(0,1)$ respectively.

With the third order RK time discretization (\ref{eq: rk3}), the last stage of (\ref{cellaverage2}) can
be written as
\beq
\bar{u}^{n+1}_h=\bar{u}^n_h-\lambda_x(\hat{H}^{rk}_{i+\f{1}{2},j}-\hat{H}^{rk}_{i-\f{1}{2},j})
-\lambda_y(\hat{G}^{rk}_{i,j+\f{1}{2}}-\hat{G}^{rk}_{i,j-\f{1}{2}}),
\label{2Deuler}
\eeq
where $\lambda_x=\Delta t/h^x$ and $\lambda_y=\Delta t/h^y$. $\hat{H}^{rk}_{i+\f{1}{2},j}$
is the integral of the numerical flux $\hat H^{rk}(x_{i+\f12}, y)$ along the cell interface 
$\{x_{i+\f12}\}\times J_j$,
which could be approximated by a numerical quadrature. At each fixed quadrature point $(x_{i+\f12}, y)$, 
$\hat H^{rk}(x_{i+\f12}, y)$ is defined the same as (\ref{rkfinal}). Similarly for $\hat G^{rk}_{i,j+\f12}$.

Let $u_m=\min_{x,y} u_0(x,y)$ and $u_M=\max_{x,y} u_0(x,y)$, numerically to preserve the cell averages within the range $[u_m,u_M]$, we are looking for the type of limiters,
\begin{eqnarray}
\tilde{H}_{i+\f12, j}&=&\theta_{i+\f12, j}(\hat H^{rk}_{i+\f12, j}-\hat h_{i+\f12, j})+\hat h_{i+\f12,j}, \label{fluxmh} \\
\tilde{G}_{i, j+\f12}&=&\theta_{i, j+\f12}(\hat G^{rk}_{i, j+\f12}-\hat g_{i, j+\f12})+\hat g_{i,j+\f12},
\label{fluxmg}
\end{eqnarray}
such that
\begin{eqnarray}
\label{MP}
u_m\le (\bar u_h)^{n}_{i, j}-\lambda_x (\tilde H_{i+\f12, j}-\tilde H_{i-\f12, j})- \lambda_y (\tilde G_{i, j+\f12}-\tilde G_{i, j-\f12}) \le u_M,
\end{eqnarray}
where $\hat h$ and $\hat g$ are first order monotone fluxes which can form a maximum principle preserving first order scheme similarly as the one dimensional case.
(\ref{fluxmh})-(\ref{MP}) form coupled inequalities for the limiting parameters $\theta_{i+\f12, j}, \theta_{i, j+\f12}$.
In each cell $K_{i,j}$, as the 1D case, the MPP flux limiters can be parametrized in the sense that we can find a group of numbers $\Lambda_{L, i, j}, \Lambda_{R, i, j},\Lambda_{D, i, j},\Lambda_{U, i, j}$, such that the numerical solutions of \eqref{2Deuler} satisfy the MPP property (\ref{MP}) with
\begin{eqnarray*}
(\theta_{i-\f12, j}, \theta_{i+\f12, j},\theta_{i, j-\f12},\theta_{i, j+\f12})\in [0, \Lambda_{L, i, j}]\times [0, \Lambda_{R, i, j}]\times [0, \Lambda_{D, i, j}] \times [0,\Lambda_{U, i, j}].
\end{eqnarray*}
For the maximum value case, let
\begin{eqnarray}
\label{1mp}
\Gamma^M_{i, j}=u_{M}-\left((\bar u_h)^n_{i, j}-\lambda_x (\hat h_{i+\f12, j}-\hat h_{i-\f12, j})- \lambda_y (\hat g_{i, j+\f12}-\hat g_{i, j-\f12})\right) \ge 0,
\end{eqnarray}
when a monotone numerical flux is used under a suitable CFL constraint, which will be specified in the numerical part. Denote
\begin{eqnarray}
\begin{cases} F_{i-\f12, j}=\lambda_x (\hat H^{rk}_{i-\f12, j} -\hat h_{i-\f12, j}),\\
F_{i+\f12, j}= -\lambda_x (\hat H^{rk}_{i+\f12, j} -\hat h_{i+\f12, j}),\\
F_{i, j-\f12}=\lambda_y (\hat G^{rk}_{i, j-\f12} -\hat g_{i, j-\f12}),\\
F_{i, j+\f12}=-\lambda_y (\hat G^{rk}_{i, j+\f12} -\hat g_{i, j+\f12}).
\end{cases}
\end{eqnarray}
The coupled inequalities (\ref{fluxmh})-(\ref{MP}) can be rewritten as
\begin{eqnarray}
\label{cMP}
\theta_{i+\f12, j}  F_{i+\f12, j}+\theta_{i-\f12, j}  F_{i-\f12, j}+\theta_{i, j+\f12}  F_{i, j+\f12}+\theta_{i, j-\f12}  F_{i, j-\f12}\le \Gamma^M_{i, j},
\end{eqnarray}
To decouple the inequality (\ref{cMP}), for the specific cell $K_{i,j}$, two steps are followed:
\begin{enumerate}
\item Identify positive values out of the four locally defined numbers $F_{i-\f12, j}$, $F_{i+\f12, j}$, $F_{i, j-\f12}$,  $F_{i, j+\f12}$;
\item Corresponding to those positive values, collectively, the limiting parameters can be defined. For example,
if $F_{i+\f12, j}, F_{i-\f12, j}>0$ and $F_{i, j-\f12}, F_{i, j+\f12}\le0$, then
\begin{eqnarray}
\label{UC}
\begin{cases} \Lambda^M_{i+\f12, j}, \Lambda^M_{i-\f12, j}=\min (\frac{\Gamma^M_{i,j}}{F_{i+\f12, j}+F_{i-\f12, j}}, 1),\\
	\Lambda^M_{i, j-\f12}, \Lambda^M_{i, j+\f12}=1.
\end{cases}
\end{eqnarray}
\end{enumerate}
For the minimum value part, let
\begin{eqnarray}
\label{ump}
\Gamma^m_{i, j}=u_{m}-\left((\bar u_h)^n_{i, j}-\lambda_x (\hat h_{i+\f12, j}-\hat h_{i-\f12, j})- \lambda_y (\hat g_{i, j+\f12}-\hat g_{i, j-\f12})\right) \le 0.
\end{eqnarray}
The coupled inequalities  (\ref{fluxmh})-(\ref{MP}) can be rewritten as
\begin{eqnarray}
\label{cuMP}
\Gamma^m_{i, j} \le \theta_{i+\f12, j}  F_{i+\f12, j}+\theta_{i-\f12, j}  F_{i-\f12, j}+\theta_{i, j+\f12}  F_{i, j+\f12}+\theta_{i, j-\f12}  F_{i, j-\f12}.
\end{eqnarray}
A similar procedure would be applied:
\begin{enumerate}
\item Identify negative values out of the four locally defined numbers $F_{i-\f12, j}$, $F_{i+\f12, j}$, $F_{i, j-\f12}$,  $F_{i, j+\f12}$;
\item Corresponding to the negative values, collectively, the limiting parameters can be defined. For example, if $F_{i, j-\f12}, F_{i, j+\f12}\ge 0$ and $F_{i-\f12, j}, F_{i+\f12, j}<0$, then
\begin{eqnarray}
\label{LC}
\begin{cases} \Lambda^m_{i-\f12, j}, \Lambda^m_{i+\f12, j}=\min(\frac{\Gamma^m_{i,j}}{F_{i-\f12, j}+F_{i+\f12, j}}, 1)\\
	\Lambda^m_{i, j-\f12}, \Lambda^m_{i, j+\f12}=1.
\end{cases}
\end{eqnarray}
\end{enumerate}
Namely, all high order fluxes which possibly contribute (beyond that of the first order fluxes) to the overshooting or undershooting of the updated value shall be limited by the same scaling. Similarly we can find $\Lambda^M_{i, j\pm\f12}$ and $\Lambda^m_{i, j\pm\f12}$, The range of the limiting parameters satisfying MPP for the cell average in cell $K_{i,j}$ therefore can be defined by
\begin{eqnarray}
\label{LocalLM}
\begin{cases}
\Lambda_{L, i, j}=\min(\Lambda^M_{i-\f12, j}, \Lambda^m_{i-\f12, j}),\\
\Lambda_{R, i, j}=\min(\Lambda^M_{i+\f12, j}, \Lambda^m_{i+\f12, j}),\\
\Lambda_{U, i, j}=\min(\Lambda^M_{i, j+\f12} , \Lambda^m_{i, j+\f12}),\\
\Lambda_{D, i, j}=\min(\Lambda^M_{i, j-\f12} , \Lambda^m_{i, j-\f12}).\\
\end{cases}
\end{eqnarray}
Considering the limiters from neighboring nodes, finally the local limiting parameters are defined to be
\begin{eqnarray}
\label{2Dlmt}
\begin{cases} \theta_{i+\f12, j}=\min(\Lambda_{R, i, j},  \Lambda_{L, i+1, j}),\\
\theta_{i, j+\f12}=\min(\Lambda_{U, i, j} , \Lambda_{D, i, j+1}).
\end{cases}
\end{eqnarray}

%% file: simulations.tex
\section{Numerical simulations}
\label{sec3}
\setcounter{equation}{0}
\setcounter{figure}{0}
\setcounter{table}{0}


In this section, we apply the parametrized MPP flux limiter to the DG method 
for solving several convection-diffusion problems and the incompressible
Navier-Stokes equations. The method is denoted as ``MPPDG'', whereas the original
DG method without the MPP flux limiter is denoted as ``DG''. The DG method is coupled with the third
order SSP RK time discretization (\ref{eq: rk3}). The time step size in this paper is defined by
\beq
\Delta t = \min\left(\f{CFLC}{\max|f'(u)|} h, \f{CFLD}{\max|a'(u)|} h^2 \right)
\label{cfl1}
\eeq
for the one-dimensional case (\ref{ad1}) and
\beq
\Delta t = \min\left(\f{CFLC}{\max|f'(u)|/h^x+ \max|g'(u)|/h^y},
\f{CFLD/\Lambda}{1/(h^x)^2+1/(h^y)^2} \right) 
\label{cfl2}
\eeq
for the two-dimensional case (\ref{ad2}), where $\Lambda$ is the maximum absolute eigenvalue of matrix $\mathbf{A}$ in (\ref{flux2d}). Here ``CFLC'' corresponds the CFL number for the convection
part, and ``CFLD'' corresponds the CFL number for the diffusion part which should be small enough. In particular, in the following, we take $CFLC=0.3, 0.18, 0.1$ from \cite{cockburn2001runge} and $CFLD=0.06, 0.01, 0.005$ as in \cite{wang2013error}, for DG method with $P^1$, $P^2$ and $P^3$ polynomial spaces respectively, unless otherwise specified. 
$\alpha$ in (\ref{flux1}) and (\ref{flux2d}) are chosen to be $1$ for $P^1$ and $10$ for $P^2$ and $P^3$. Each problem is computed to the final time ``T'' on the mesh of ``$N$'' cells for the one-dimensional case and ``$N^2$'' cells for the two dimensional case. For solutions with discontinuity, the TVB limiter \cite{cockburn2001runge} with a parameter $M_{tvb}$ usually needs to be applied to ensure stability. For some of the following cases, we avoid the TVB limiter to see the good performance of the MPP flux limiter, if the numerical solutions are still stable without the TVB limiter. For all figures, the cell averages of the numerical solutions are displayed.

\subsection{Basic tests of MPP for the one dimensional case}
\label{sec4.1}
\begin{exa}(Accuracy test)
We first test the accuracy for the linear equation
\beq
u_t+u_x=\varepsilon u_{xx},
\label{eq401}
\eeq
with initial condition $u(x,0)=\sin^4(x)$ on $[0, 2\pi]$ and periodic boundary conditions. 
The exact solution is
\beq
u(x,t)=\f{3}{8}-\f{1}{2}\exp(-4 \varepsilon t)\cos(2(x-t))+\f{1}{8}\exp(-16 \varepsilon t)\cos(4(x-t)).
\label{eq402}
\eeq
Let $\varepsilon=0.0001$ and the final time $T=1$, we show the $L^1$ and $L^\infty$ errors and orders for $P^2$ and $P^3$ cases in Table \ref{tab401}. For the $P^2$ case, the MPP flux limiter can limit the undershoot (negative minimum values without the MPP limiter) within the theoretical bounds, without affecting the overall accuracy, since clear 3rd order accuracy for both DG and MPPDG are observed. For $P^3$ case, there is no overshoot or undershoot of the DG solution, thus the flux limiter are not effective.

\begin{table}[ht]
\centering
\caption{$L^1$ and $L^\infty$ errors and orders for (\ref{eq401}) with initial 
condition $u(x,0)=\sin^4(x)$ and exact solution (\ref{eq402}). $T=1$. The time step
is (\ref{cfl1}) for $P^2$ and it is $\Delta t = \min\left(\f{CFLC}{\max|f'(u)|} h^{4/3}, \f{CFLD}{\max|a'(u)|} h^2 \right)$ for $P^3$ here.}
\vspace{0.2cm}
  \begin{tabular}{|c||c|c|c|c|c|c|c|}
    \hline
 &    N &  $L^1$ error  & order  & $L^\infty$ error & order  & $(\bar u_h)_{min}$    & $(\bar u_h)_{max}$   \\ \hline
\multirow{5}{*}{$P^2$ DG}
 &   16 &     1.61E-03 &       --&     6.25E-03 &       --& -0.0004060923125&  0.9727611964447 \\ \cline{2-8}
 &   32 &     1.87E-04 &     3.10&     8.19E-04 &     2.93&  0.0001429702538&  0.9806422360262 \\ \cline{2-8}
 &   64 &     2.30E-05 &     3.03&     1.04E-04 &     2.97&  0.0000033153325&  0.9960977566409 \\ \cline{2-8}
 &  128 &     2.86E-06 &     3.01&     1.30E-05 &     3.00& -0.0000000980932&  0.9991190862113 \\ \cline{2-8}
 &  256 &     3.59E-07 &     2.99&     1.61E-06 &     3.02&  0.0000001699923&  0.9994283533012 \\ \hline
\multirow{5}{*}{$P^2$ MPPDG}
 &   16 &     1.56E-03 &       --&     6.25E-03 &       --&  0.0000000000000&  0.9727609449833 \\ \cline{2-8}
 &   32 &     1.86E-04 &     3.07&     8.19E-04 &     2.93&  0.0001617902828&  0.9806422317515 \\ \cline{2-8}
 &   64 &     2.29E-05 &     3.02&     1.04E-04 &     2.97&  0.0000062375238&  0.9960977566566 \\ \cline{2-8}
 &  128 &     2.86E-06 &     3.00&     1.30E-05 &     3.00&  0.0000000745687&  0.9991190862113 \\ \cline{2-8}
 &  256 &     3.59E-07 &     2.99&     1.61E-06 &     3.02&  0.0000001699923&  0.9994283533012 \\ \hline
\multirow{5}{*}{$P^3$ DG}
 &   16 &     1.26E-04 &       --&     4.22E-04 &       --&  0.0003066194525&  0.9737481094672 \\ \cline{2-8}
 &   32 &     8.13E-06 &     3.95&     2.64E-05 &     4.00&  0.0001895991922&  0.9807169190008 \\ \cline{2-8}
 &   64 &     5.03E-07 &     4.01&     1.73E-06 &     3.94&  0.0000087482084&  0.9961066119535 \\ \cline{2-8}
 &  128 &     3.11E-08 &     4.02&     1.07E-07 &     4.01&  0.0000005241038&  0.9991201152638 \\ \cline{2-8}
 &  256 &     1.90E-09 &     4.03&     6.35E-09 &     4.07&  0.0000002391899&  0.9994284684069 \\ \hline
\multirow{5}{*}{$P^3$ MPPDG}
 &   16 &     1.26E-04 &       --&     4.22E-04 &       --&  0.0003066194525&  0.9737481094672 \\ \cline{2-8}
 &   32 &     8.13E-06 &     3.95&     2.64E-05 &     4.00&  0.0001895991922&  0.9807169190008 \\ \cline{2-8}
 &   64 &     5.03E-07 &     4.01&     1.73E-06 &     3.94&  0.0000087482084&  0.9961066119535 \\ \cline{2-8}
 &  128 &     3.11E-08 &     4.02&     1.07E-07 &     4.01&  0.0000005241038&  0.9991201152638 \\ \cline{2-8}
 &  256 &     1.90E-09 &     4.03&     6.35E-09 &     4.07&  0.0000002391899&  0.9994284684069 \\ \hline
  \end{tabular}
\label{tab401}
\end{table}
\end{exa}

\begin{exa}
In the second example, we consider a linear advection equation \cite{Jiang_Shu}
\beq
u_t+u_x=0, \quad u(x,0)=u_0(x)
\label{eq420}
\eeq
with 
\beqa
u_0(x)=
\begin{cases}
\frac16(G(x,\beta,z-\delta)+G(x,\beta,z+\delta)+4G(x,\beta,z)), &\quad -0.8\le x\le -0.6; \\
1, &\quad -0.4\le x \le -0.2; \\
1-|10(x-0.1)|, &\quad 0\le x \le 0.2; \\
\frac16(F(x,\gamma,a-\delta)+F(x,\gamma,a+\delta)+4F(x,\gamma,a)), &\quad 0.4\le x \le 0.6; \\
0, &\quad \text{ otherwise.}
\end{cases} \label{eq421}\\
G(x,\beta,z)=e^{-\beta(x-z)^2}, \qquad F(x,\gamma,a)=\sqrt{\max(1-\gamma^2(x-a)^2,0)} \nonumber
\eeqa
where the constants are taken as $a=0.5$, $z=-0.7$, $\delta=0.005$, $\gamma=10$ and $\beta=\log2/36\delta^2$.
The computational domain is $[-1, 1]$ with periodic boundary condition. The solution contains a smooth 
but narrow combinations of Gaussians, a square wave, a sharp triangle wave and a half ellipse. In Fig. \ref{fig410a}, we show the $P^2$ solution at $T=8$ with mesh $N=200$. The TVB limiter with $M_{tvb}=10$ is used. The minimum and maximum values without the MPP flux limiter are $-0.0000872949879$ and $1.0000844689587$, while with the MPP flux limiter they are $0.0000000000009$ and $0.9999992542962$.
In Fig. \ref{fig410b} and the zoom-in Figs. \ref{fig410c} and \ref{fig410d}, we show the $P^2$ solution without the TVB limiter, in which the effect of the MPP flux limiter can be clearly observed. 

\begin{figure}[htbp]
\begin{center}
\subfigure[with TVB limiter]{
\includegraphics[scale=0.3]{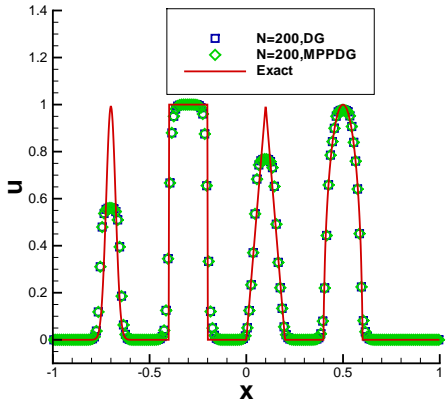}\label{fig410a}},
\subfigure[without TVB limiter]{
\includegraphics[scale=0.3]{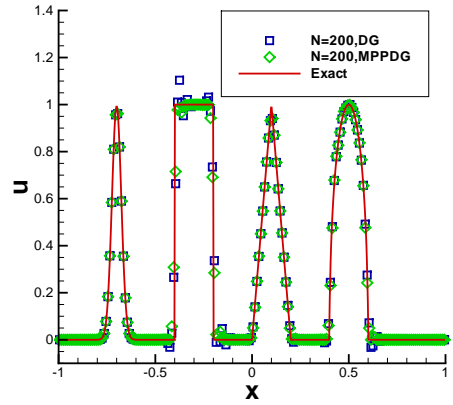}\label{fig410b}} \\
\subfigure[zoom-in of (b) at the bottom]{
\includegraphics[scale=0.3]{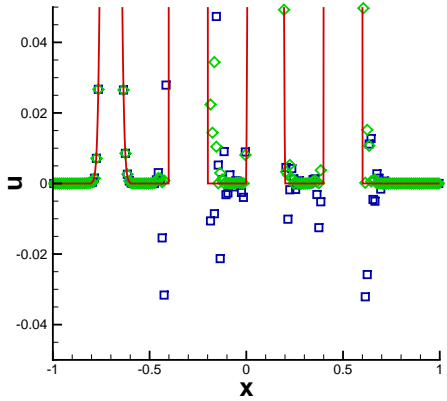}\label{fig410c}},
\subfigure[zoom-in of (b) on the top]{
\includegraphics[scale=0.3]{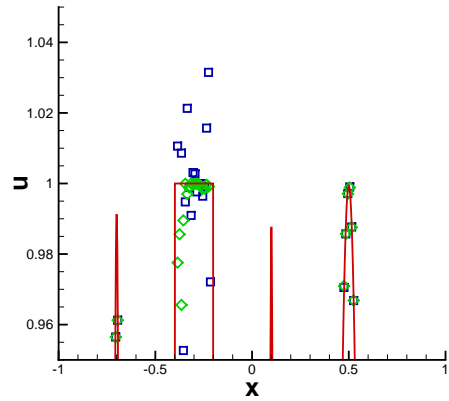}\label{fig410d}}
\end{center}
\caption{Linear advection equation (\ref{eq420}) with initial condition (\ref{eq421}). $T=8$. 
Solid line: the exact solution; Symbols: cell averages of the $P^2$ numerical solutions with mesh $N=200$.}
\label{fig410}
\end{figure} 
\end{exa}

\begin{exa} (Porous medium equation)
This is a typical example of the degenerate parabolic equations \cite{muskat1946flow}. We consider
\beq
u_t=(u^m)_{xx}, \qquad m>1.
\label{eq403}
\eeq
The Barenblatt solution
\beq
B_m(x,t)=t^{-s}\left[\left(1-\f{s(m-1)}{2m}\f{|x|^2}{t^{2s}}\right)_+\right]^{1/(m-1)}
\label{eq404}
\eeq
is an exact solution to (\ref{eq403}) with compact support, where $v_+=\max(v,0)$ and $s=(m+1)^{-1}$.
The initial condition is $B_m(x,1)$. We compute the numerical solutions with $m=2, 3, 5, 8$ to the time $T=2$ with zero boundary conditions on $[-6, 6]$. For this example, to see the difference between DG and MPPDG methods, we use the $P^3$ piecewise polynomial space and the TVB limiter with $M_{tvb}=1$. With $N=80$, in Table \ref{tab402}, we can clearly observe the negative undershoots for the DG solutions. There is no such negative undershoot in the MPPDG solutions. The corresponding MPPDG numerical solutions are plotted in Fig. \ref{fig401}, which match the Barenblatt solution very well. 

\begin{figure}[htbp]
\begin{center}
\subfigure[m=2]
{\includegraphics[scale=0.3]{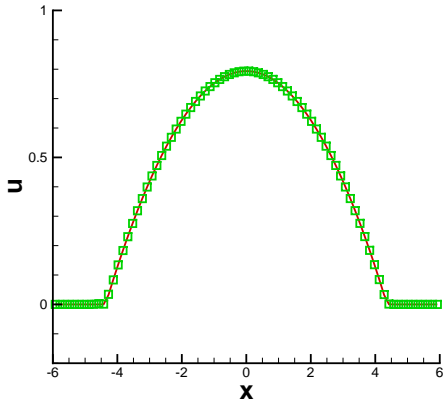}},
\subfigure[m=3]
{\includegraphics[scale=0.3]{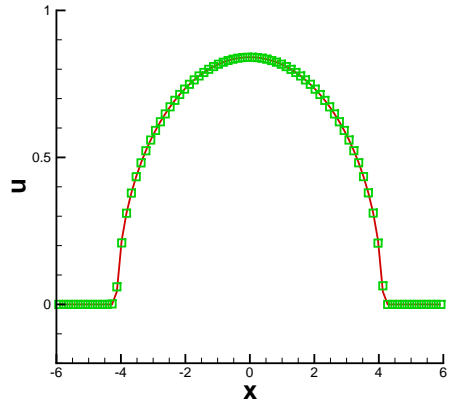}} \\
\subfigure[m=5]
{\includegraphics[scale=0.3]{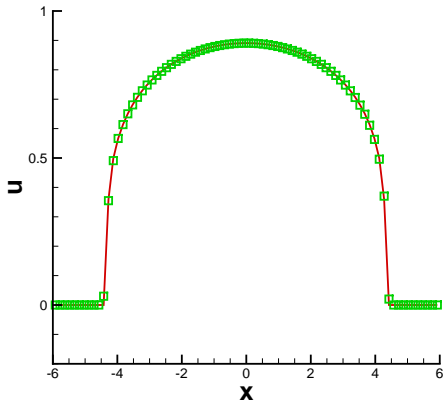}},
\subfigure[m=8]
{\includegraphics[scale=0.3]{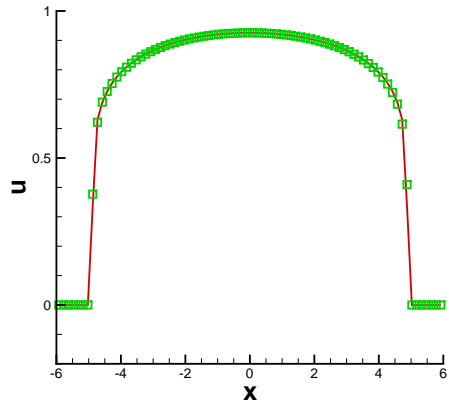}}
\end{center}
\caption{Porous medium equation (\ref{eq403}) with the Barenblatt solution (\ref{eq404}) at $T=2$. The solid line denotes the Barenblatt solution and the symbols are $P^3$ MPPDG numerical solutions, $N=80$.}
\label{fig401}
\end{figure}

\begin{table}[ht]
\centering
\caption{Minimum values of the $P^3$ solutions for Porous medium equation (\ref{eq403}) at $T=2$, $N=80$.}
\vspace{0.2cm}
  \begin{tabular}{|c|c|c|}
    \hline
    m & $(\bar u_h)_{min}$ of DG  & $(\bar u_h)_{min}$ of MPPDG   \\ \hline
    2 &    -0.0000158453675  &     0.0000000000000 \\ \hline
    3 &    -0.0000069345796  &     0.0000000000000 \\ \hline
    5 &    -0.0000000026392  &     0.0000000000000 \\ \hline
    8 &     0.0000000000000  &     0.0000000000000 \\ \hline
  \end{tabular}
\label{tab402}
\end{table}
\end{exa}

\begin{exa} (Buckley-Leverett equation)
Now we consider the Buckley-Leverett convection-diffusion equation, which is a model often
used in reservoir simulations \cite{leveque2002finite}
\beq
u_t+f(u)_x=\varepsilon (\nu(u)u_x)_x.
\label{eq405}
\eeq
We take $\varepsilon=0.01$ and boundary conditions $u(0,t)=1$ and $u(1,t)=0$ on $[0, 1]$. The function $\nu(u)$ and the initial condition are given as
\beq
\nu(u)=
\begin{cases}
4u(1-u),& 0\le{u}\le{1}, \\
0, &\text{otherwise},
\end{cases}
\qquad
u(x,0)=
\begin{cases}
1-3x, &\quad 0\le{x}\le{\f{1}{3}}, \\
0, &\quad \f{1}{3}\le{x}\le{1},
\end{cases}
\eeq
with an $s$-shape function 
\[
f(u)=\f{u^2}{u^2+(1-u)^2}.
\]
In Fig. \ref{fig402}, we show the numerical solutions of DG and MPPDG methods at $T=0.2$ on the mesh of $N=100$ compared with the reference solution of MPPDG on the mesh of $N=500$ for $P^k$, $k=1, 2, 3$ respectively. We use the TVB limiter with $M_{tvb}=10$. All solutions match each other well. However, the DG method would have negative undershoots
while MPPDG does not, which can be seen from Table \ref{tab403} and the zoom-in figure in Fig. \ref{fig402} (b) for the $P^1$ case.

\begin{figure}[htbp]
\begin{center}
\subfigure[$P^1$]
{\includegraphics[scale=0.3]{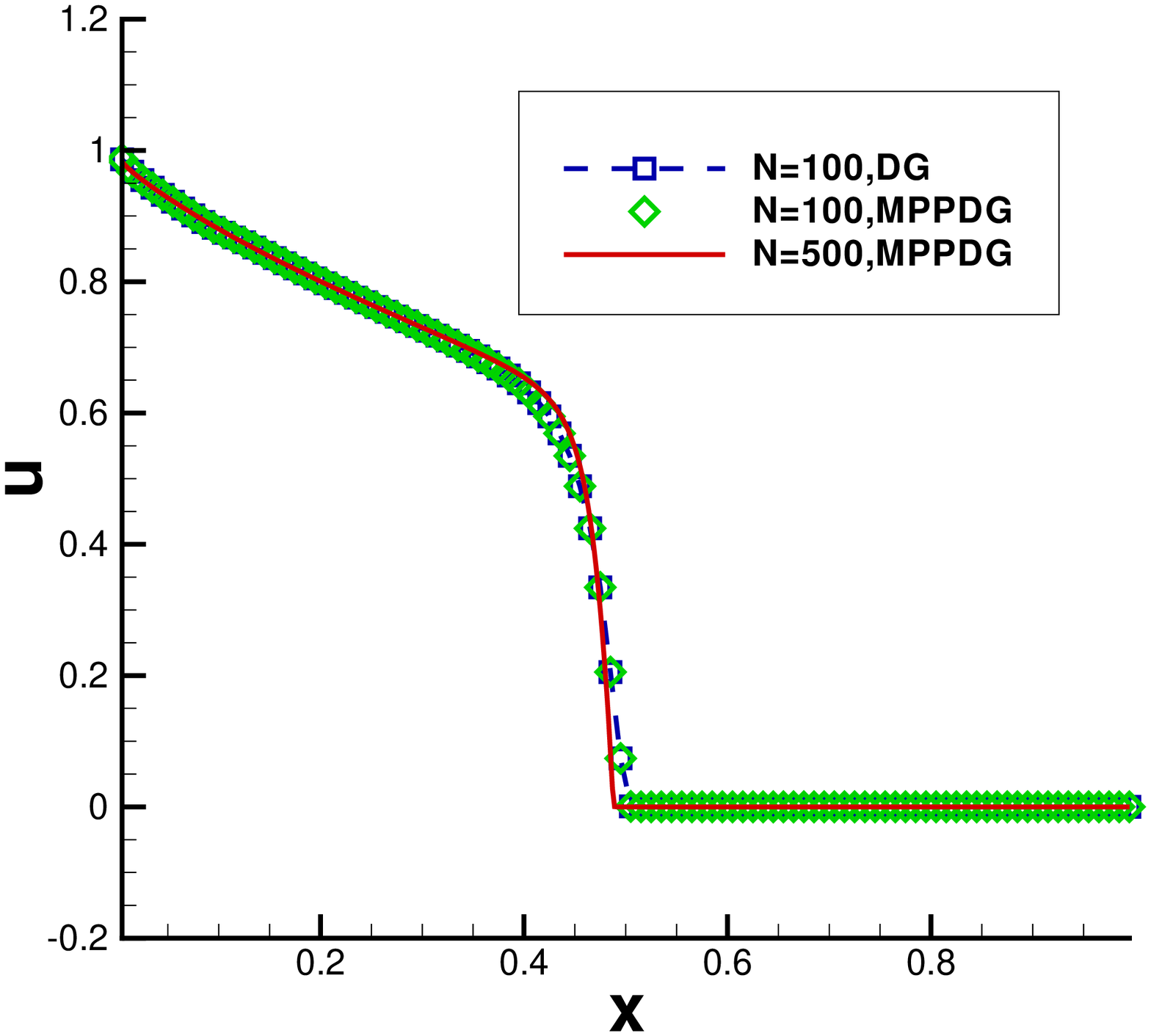}},
\subfigure[Zoom in near the undershoot of (a)]
{\includegraphics[scale=0.3]{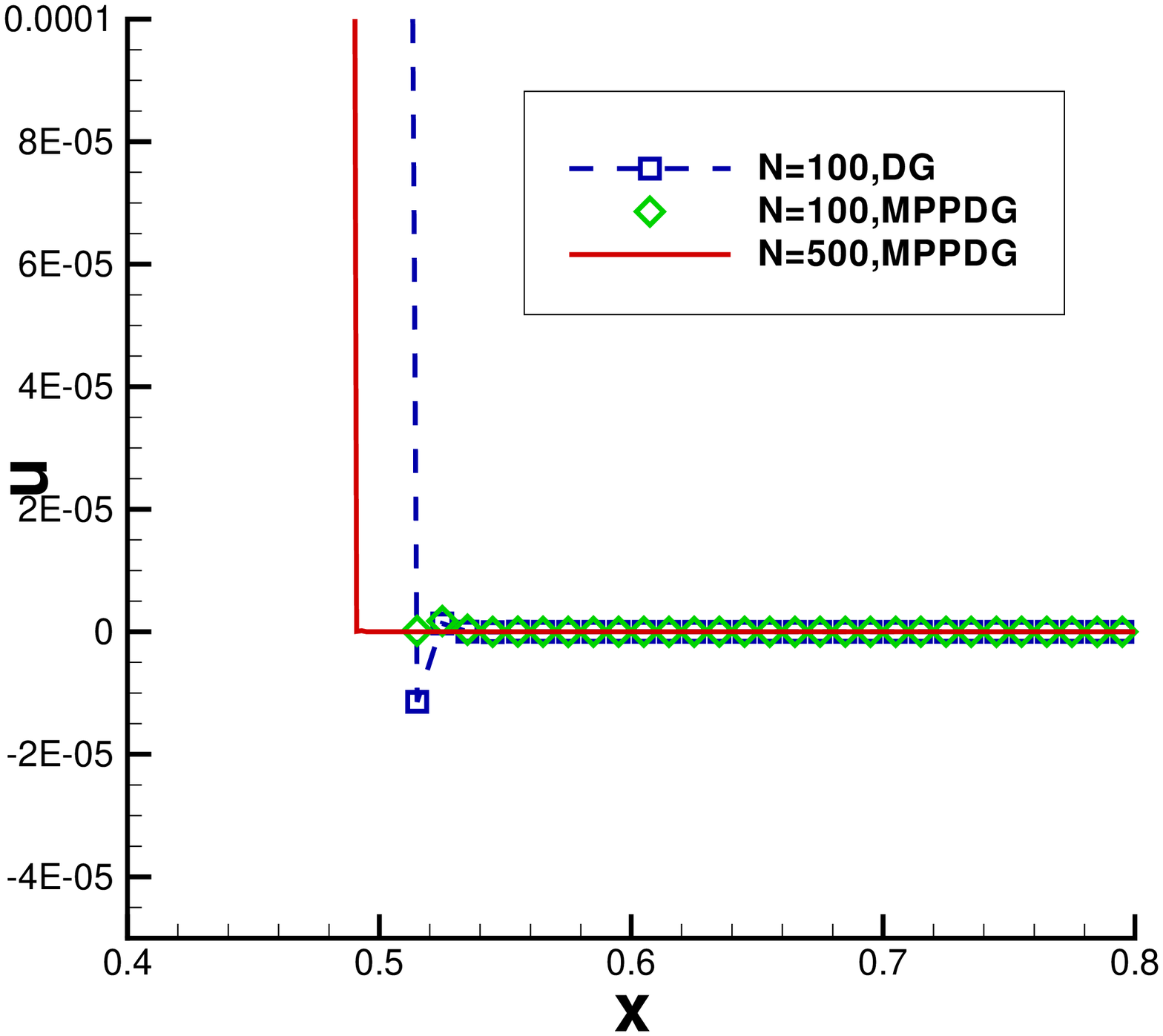}}\\
\subfigure[$P^2$]
{\includegraphics[scale=0.3]{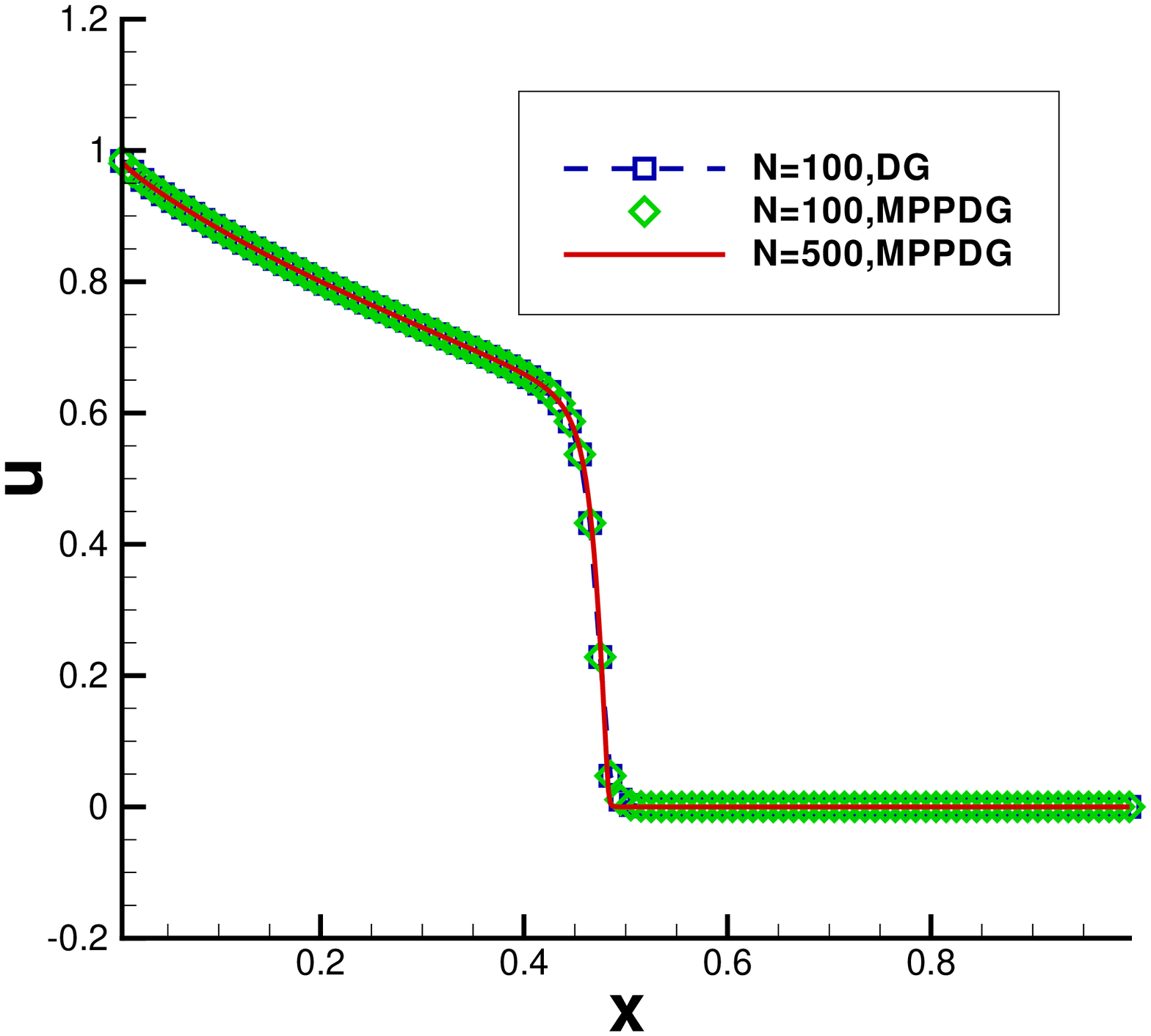}},
\subfigure[$P^3$]
{\includegraphics[scale=0.3]{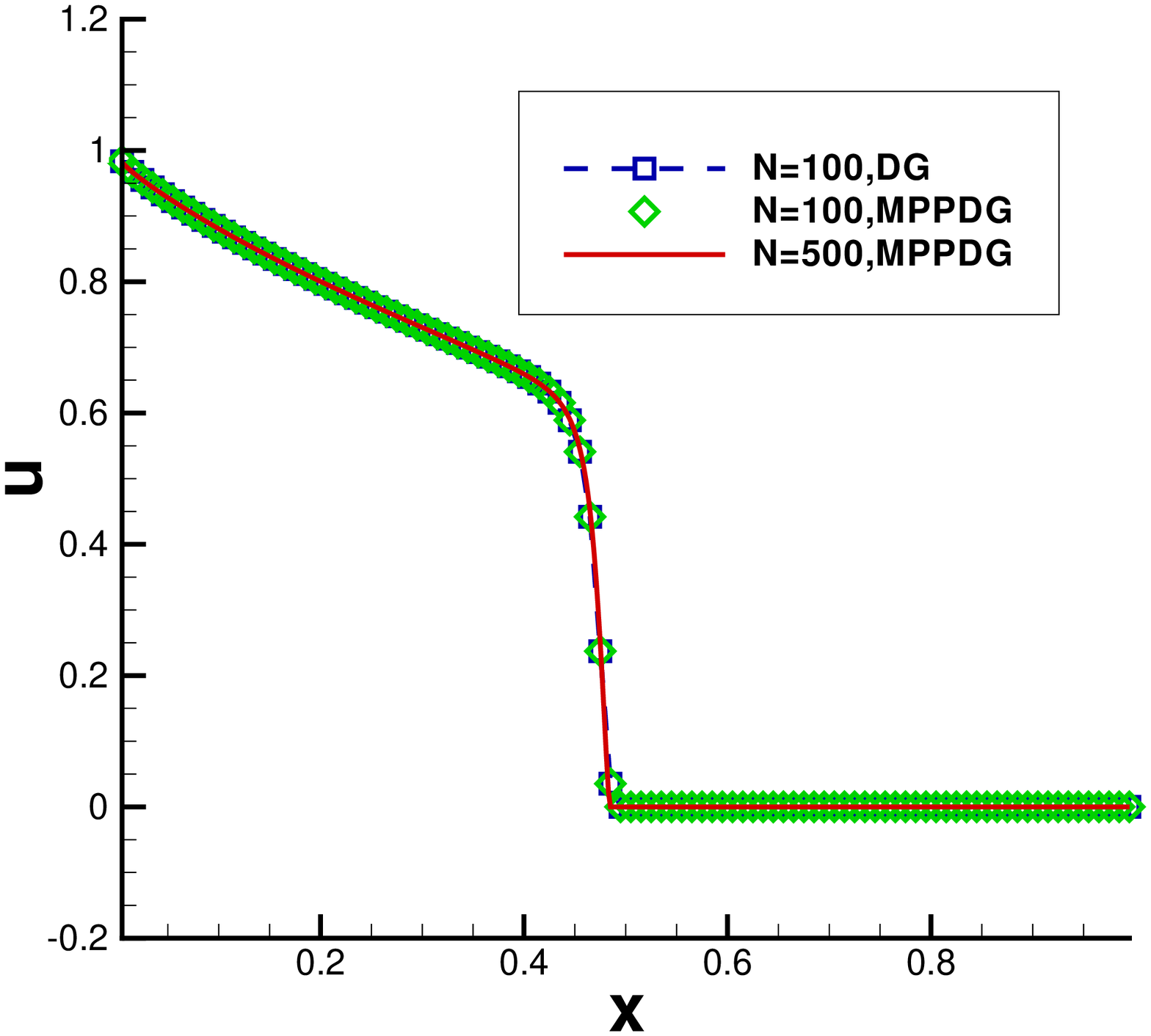}}
\end{center}
\caption{Numerical solutions of Buckley-Leverett equation (\ref{eq405}) at $T=0.2$. }
\label{fig402}
\end{figure}

\begin{table}[ht]
\centering
\caption{Minimum values of $P^k$ solutions for Buckley-Leverett equation (\ref{eq405}) at $T=0.2$, $N=100$.}
\vspace{0.2cm}
  \begin{tabular}{|c|c|c|}
  \hline
 $P^k$ & $(\bar u_h)_{min}$ of DG  & $(\bar u_h)_{min}$ of MPPDG   \\ \hline
 $k=1$ &    -0.0000114304519  &     0.0000000000000 \\ \hline
 $k=2$ &     0.0000000000000  &     0.0000000000000 \\ \hline
 $k=3$ &    -0.0000000462710  &     0.0000000000000 \\ \hline
  \end{tabular}
\label{tab403}
\end{table}
\end{exa}

\subsection{Basic tests of MPP for two dimensional case}
\label{sec4.2}
\begin{exa}(Accuracy test)
We test the linear equation similarly as the one dimensional case
\beq
u_t+u_x+u_y=\varepsilon(u_{xx}+u_{yy}), \quad u(x,y,0)=\sin^4(x+y),
\label{eq406}
\eeq
on $[0, 2\pi]^2$ with periodic boundary conditions. The
exact solution is
\beq
u(x,y,t)=\f{3}{8}-\f{1}{2}\exp(-8\varepsilon t)\cos(2(x+y-2t))+
\f{1}{8}\exp(-32\varepsilon t)\cos(4(x+y-2t)).
\label{eq407}
\eeq
We take $\varepsilon=0.0001$ and $T=0.5$. We show the $L^1$ and $L^\infty$ errors and orders for $P^2$ case in Table \ref{tab404}. As expected, 3rd order accuracies have been observed for both DG and MPPDG solutions. The MPP flux limiter can limit the undershoot within the theoretical bounds, without affecting the overall accuracy.

\begin{table}[ht]
\centering
\caption{$L^1$ and $L^\infty$ errors and orders for (\ref{eq406}) with exact solution (\ref{eq407}). $P^2$ and $T=0.5$.}
\vspace{0.2cm}
  \begin{tabular}{|c||c|c|c|c|c|c|c|}
    \hline
 &   $N^2$ &  $L^1$ error  & order  & $L^\infty$ error & order  & $(\bar u_h)_{min}$    & $(\bar u_h)_{max}$   \\ \hline
\multirow{5}{*}{DG}
 & $8^2$ &     5.17E-02 &       --&     1.88E-01 &       --&  0.0280510290327&  0.7430134210472 \\ \cline{2-8}
 &$16^2$ &     6.80E-03 &     2.93&     5.29E-02 &     1.83&  0.0049275726860&  0.8938493550185 \\ \cline{2-8}
 &$32^2$ &     7.22E-04 &     3.24&     8.73E-03 &     2.60& -0.0001491953901&  0.9859512188155 \\ \cline{2-8}
 &$64^2$ &     8.55E-05 &     3.08&     1.12E-03 &     2.97& -0.0000014782832&  0.9957228872298 \\ \cline{2-8}
 &$128^2$&     1.05E-05 &     3.02&     1.41E-04 &     2.99&  0.0000014191786&  0.9981330285231 \\ \hline
\multirow{5}{*}{MPPDG}
 & $8^2$ &     5.17E-02 &       --&     1.88E-01 &       --&  0.0280510290327&  0.7430134210472 \\ \cline{2-8}
 &$16^2$ &     6.56E-03 &     2.98&     5.29E-02 &     1.83&  0.0053726995812&  0.8938465321283 \\ \cline{2-8}
 &$32^2$ &     7.17E-04 &     3.19&     8.73E-03 &     2.60&  0.0000000000001&  0.9859512187802 \\ \cline{2-8}
 &$64^2$ &     8.53E-05 &     3.07&     1.12E-03 &     2.97&  0.0000019904149&  0.9957228872298 \\ \cline{2-8}
 &$128^2$&     1.05E-05 &     3.02&     1.41E-04 &     2.99&  0.0000014191786&  0.9981330285231 \\ \hline
  \end{tabular}
\label{tab404}
\end{table}
\end{exa}

\begin{exa}(Porous medium equation)
The two-dimensional porous medium equation is
\beq
u_t=(u^2)_{xx}+(u^2)_{yy},
\label{eq408}
\eeq
with the initial condition
\beq
u(x,y,0)=
\begin{cases}
1, \quad & \text{ if } (x,y)\in [-\f{1}{2},\f{1}{2}]\times[-\f{1}{2},\f{1}{2}], \\
0, \quad & \text{ otherwise}.
\end{cases}
\label{eq409}
\eeq
on the domain $[0,1]^2$ and periodic boundary conditions. With $P^1$ piecewise polynomial space, we compare the results for DG and MPPDG methods in Figure \ref{fig403}. From the zoom-in figure,
we can clearly see that the negative value for the DG method has been eliminated by the MPPDG method.
The TVB limiter with $M_{tvb}=50$ has been used. $P^2$ and $P^3$ solutions are omitted here due to similarity. 

\begin{figure}[htbp]
\begin{center}
\subfigure[Surface of $P^1$ MPPDG solution.]
{\includegraphics[scale=0.3, angle=270]{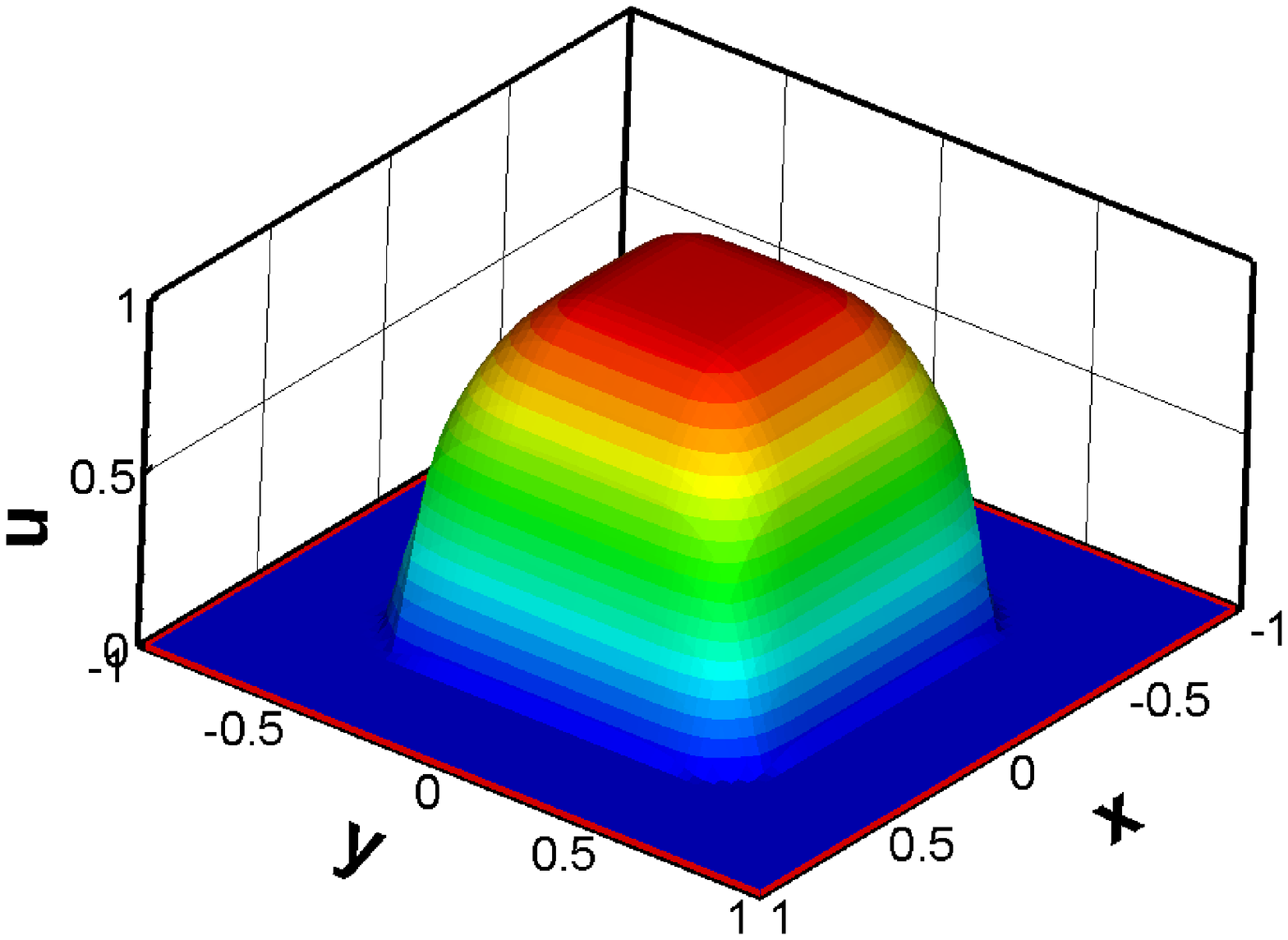}}
\subfigure[Zoom-in of cuts along $y=0$. Symbol with dashed line: DG; Solid line: MPPDG.]
{\includegraphics[scale=0.3, angle=270]{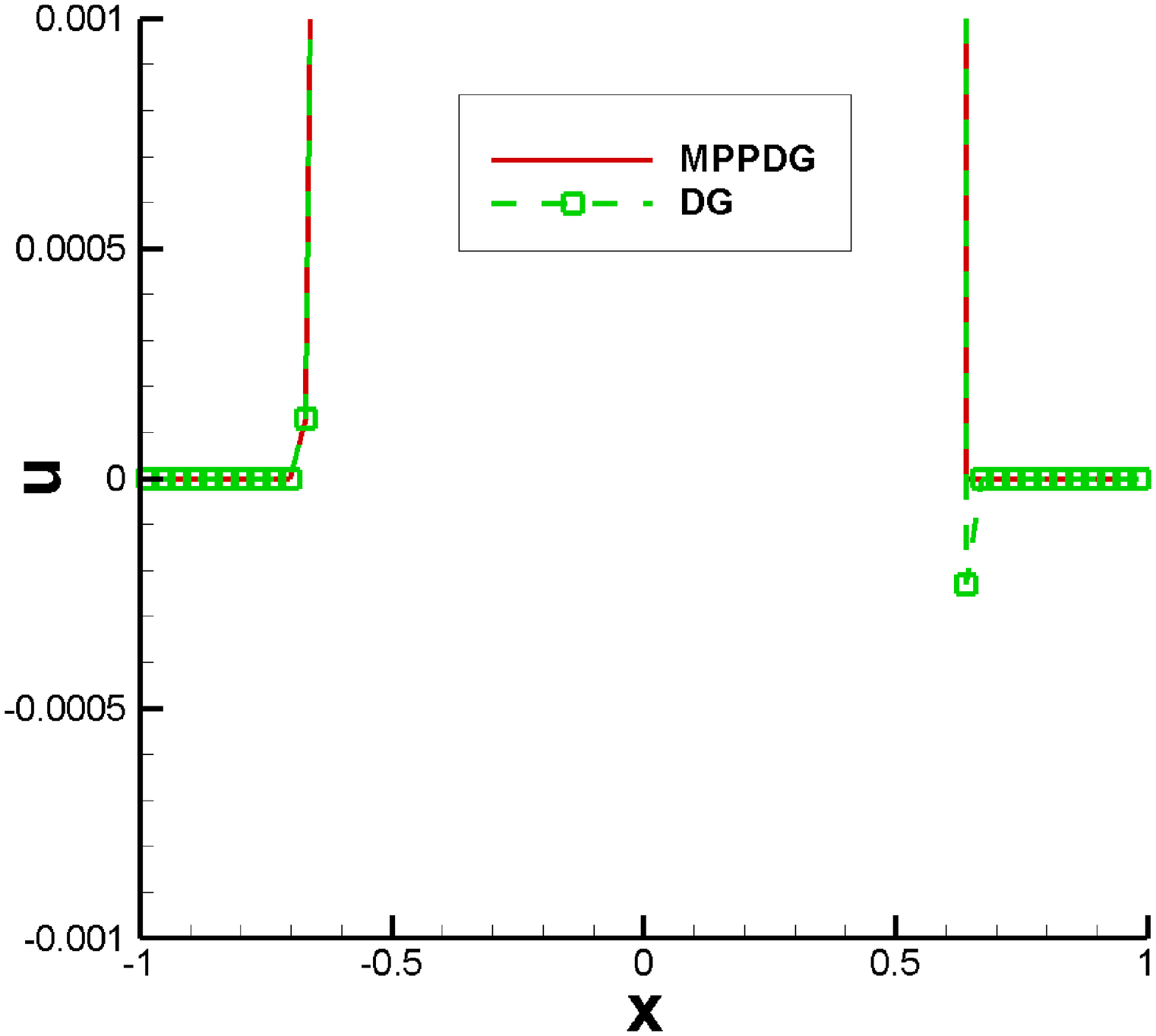}}
\end{center}
\caption{Porous medium equation (\ref{eq408}) with initial condition (\ref{eq409}). $T=0.005$. 
$P^1$ with mesh $N^2=64^2$.}
\label{fig403}
\end{figure}


\end{exa}

\begin{exa}(Buckley-Leverett equation)
The two-dimensional Buckley-Leverett equation with gravity in $y$-direction is given by \cite{jiang2013parametrized,kurganov2000new}
\beq
u_t+f(u)_x+g(u)_y=\varepsilon(u_{xx}+u_{yy}),
\label{eq410}
\eeq
where
\[
f(u)=\frac{u^2}{u^2+(1-u)^2}, \qquad g(u)=f(u)(1-5(1-u)^2).
\]
The initial condition is
\beq
u(x,y,0)=
\begin{cases}
1, & x^2+y^2<0.5, \\
0, & \text{ otherwise}.
\end{cases}
\label{eq411}
\eeq
We take $\varepsilon=0.01$ and periodic boundary conditions. We run the numerical solution to $T=0.5$ and  show the minimum and maximum values on different meshes in Table \ref{tab405}. Similarly, the TVB limiter with $M_{tvb=50}$ has been used. We can clearly see the overshoots and undershoots have been eliminated by the MPPDG method. The surface and contour for the $P^2$ MPPDG solutions on the mesh of $256^2$ grid points are displayed in Fig. \ref{fig405}. The DG solutions are also omitted here due to similarity. 

\begin{table}[ht]
\centering
\caption{Minimum and maximum values of the $P^2$ solutions for Buckley-Leverett equation (\ref{eq410}) at $T=0.5$.}
\vspace{0.2cm}
\begin{tabular}{|c||c c|c c|}
\hline
        & \multicolumn{2}{c|}{DG}                &   \multicolumn{2}{c|}{MPPDG}       \\ \hline
  $N^2$ &  $(\bar u_h)_{min}$    &  $(\bar u_h)_{max}$    &  $(\bar u_h)_{min}$   & $(\bar u_h)_{max}$    \\ \hline
 $16^2$ & -0.1692411477038  &  1.1778306576225  &  0.0000000000000 &  1.0000000000000 \\ \hline
 $32^2$ & -0.0855426715979  &  1.0518223688230  &  0.0000000000000 &  1.0000000000000 \\ \hline
 $64^2$ & -0.0317585482621  &  1.0183785242840  &  0.0000000000000 &  0.9999165339925 \\ \hline
$128^2$ & -0.0084299844740  &  1.0015417862032  &  0.0000000000000 &  1.0000000000000 \\ \hline
$256^2$ & -0.0009240813456  &  0.9999747335598  &  0.0000000000000 &  0.9999745692661 \\ \hline
\end{tabular}
\label{tab405}
\end{table}

\begin{figure}[htbp]
\begin{center}
\includegraphics[scale=0.3, angle=270]{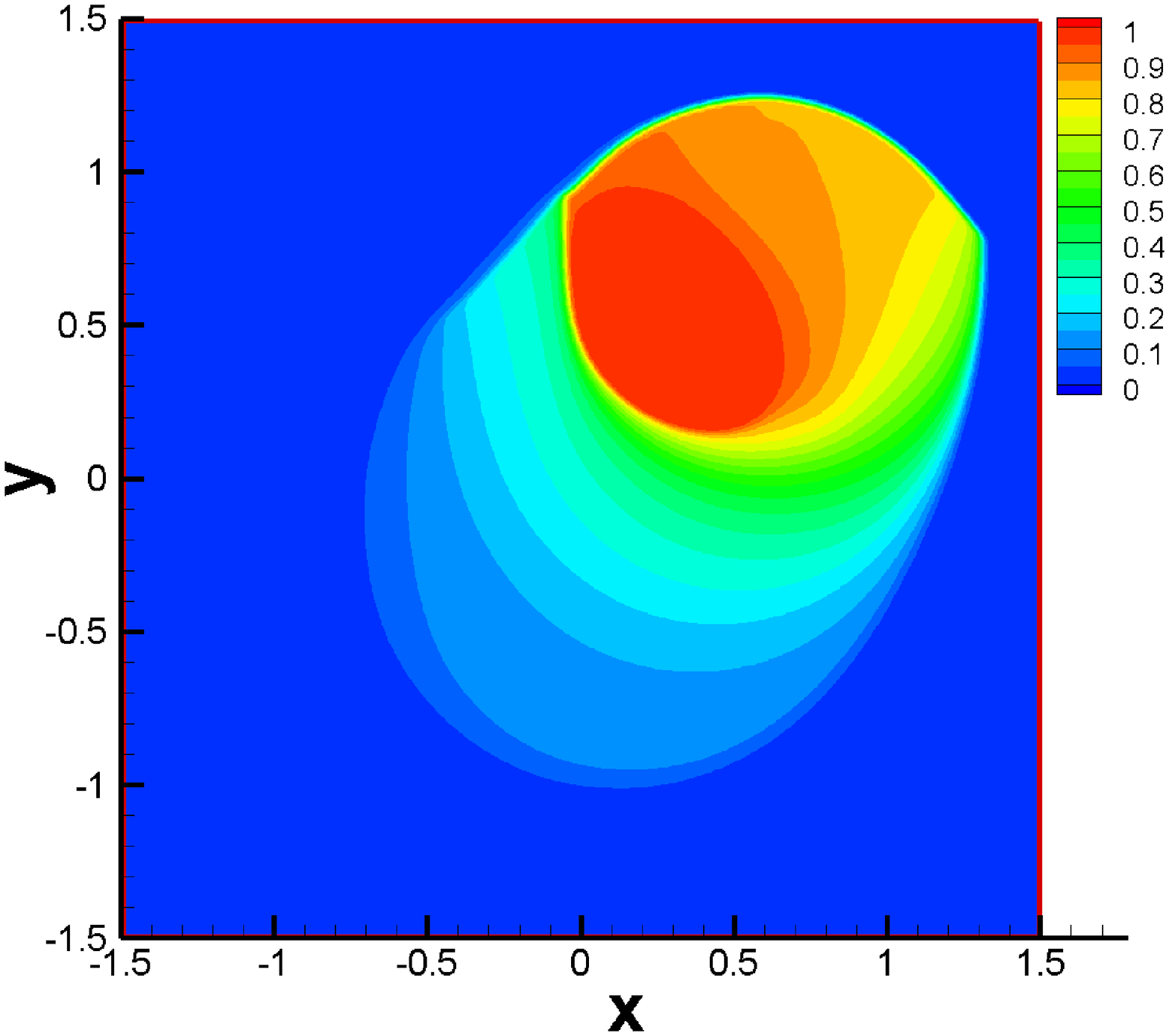}
\includegraphics[scale=0.3, angle=270]{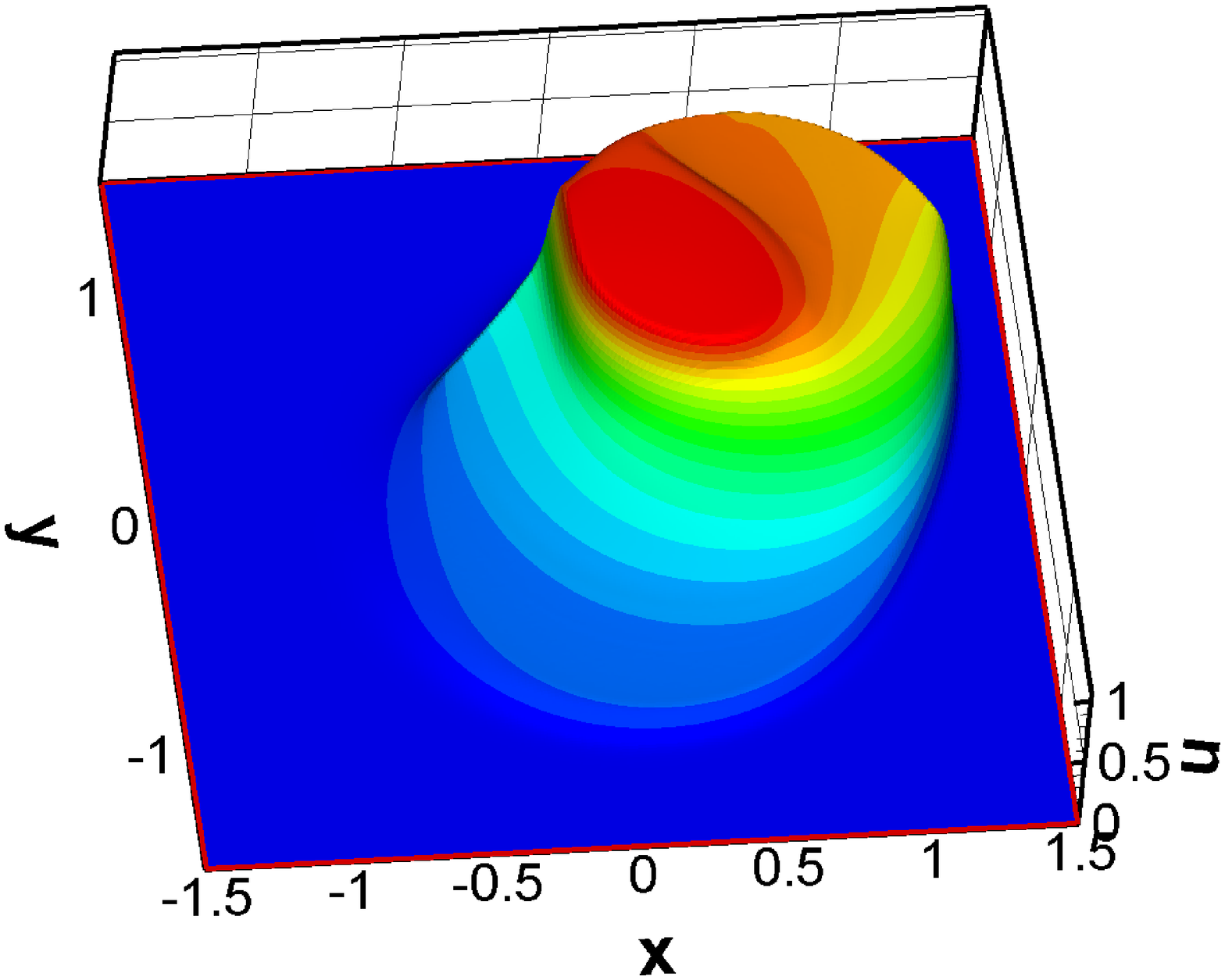}
\end{center}
\caption{Buckley-Leverett equation (\ref{eq410}) with initial condition (\ref{eq411}). $T=0.5$. 
$P^2$ MPPDG with mesh $N^2=256^2$.}
\label{fig405}
\end{figure}

\end{exa}

\subsection{Incompressible flow}
\label{sec4.3}
In this section, we consider the incompressible 
Navier-Stokes equations in the vorticity-stream function formulation
\begin{align}
\label{ins1}
\omega_t+(u\omega)_x+(v\omega)_y=\frac{1}{Re}(\omega_{xx}+\omega_{yy}), \\
\label{ins2}
\Delta \psi=\omega, \quad \langle u,v \rangle = \langle -\psi_y,\psi_x \rangle, \\
\label{ins3}
\omega(x,y,0)=\omega_0(x,y), \quad \langle u, v \rangle\cdot\mathbf{n}\text{=given on } \partial\Omega.
\end{align}
The solution to the incompressible flow problem satisfies the maximum principle
due to the divergence-free property of the velocity field. Numerically,
the discretized divergence-free condition has been delicately built into the discretization of the convection term to ensure
the MPP property of numerical solutions, see \cite{mpp_xqx} for the incompressible Euler problem. 
In the following examples, without specifying, we take $Re=100$ or $Re=\infty$ for inviscid case.
And if the TVB limiter is used, we take $M_{tvb}=50$.

\begin{exa}(Rigid body Rotation) 
We first consider an incompressible flow problem with explicitly given velocity field, 
which involves a rigid body rotation
\beq
\omega_t-(y \omega)_x+(x\omega)_y=\frac{1}{Re}(\omega_{xx}+\omega_{yy}),
\label{eq412}
\eeq
with zero boundary conditions on the domain $[-\pi, \pi]^2$. The initial condition includes
a slotted disk, a cone and a smooth hump as shown in Fig. \ref{fig406}. 
For this problem, the initial condition rotates counterclockwise. After a period of $T=2\pi$, the solution will get back to its initial position. We first take $1/Re=0$, that is without viscosity. In Fig. \ref{fig407} (left), we show the cuts along $x=0$, $y=0.8$ and $y=-2$ for the $P^2$ numerical solutions at $T=2\pi$ without the TVB limiter, we can clearly see the overshoots and undershoots have been eliminated by the MPP flux limiter. Then we take $1/Re=0.01$. In Table \ref{tab406}, the minimum and maximum values of the $P^2$ numerical solutions on different meshes at $T=0.1$ with the TVB limiter indicate that the undershoots and overshoots of the DG method can be eliminated by the MPPDG method. 


\begin{table}[ht]
\centering
\caption{Minimum and maximum values of the $P^2$ solutions for the rigid body rotation problem (\ref{eq412}) with initial condition in Fig. \ref{fig406}. $T=0.1$.}
\vspace{0.2cm}
  \begin{tabular}{|c||c c|c c|}
    \hline
        & \multicolumn{2}{c|}{DG}               &   \multicolumn{2}{c|}{MPPDG}        \\ \hline
  $N^2$ &  $(\bar \omega_h)_{min}$    &  $(\bar \omega_h)_{max}$    &  $(\bar \omega_h)_{min}$   &  $(\bar \omega_h)_{max}$   \\ \hline
  $8^2$ & -0.0086609245500  &  0.9287124547518  &  0.0000000000000 &  0.9283345801837 \\ \hline
 $16^2$ & -0.0223925807132  &  1.0211949713819  &  0.0000000000000 &  0.9999999999999 \\ \hline
 $32^2$ & -0.0250454358990  &  1.0216925918827  &  0.0000000000000 &  0.9999999999999 \\ \hline
 $64^2$ & -0.0051503335802  &  1.0054007922926  &  0.0000000000000 &  0.9999999999999 \\ \hline
$128^2$ & -0.0004058734427  &  1.0006938627099  &  0.0000000000000 &  0.9999999999999 \\ \hline
  \end{tabular}
\label{tab406}
\end{table}

\begin{figure}[htbp]
\begin{center}
\includegraphics[scale=0.3, angle=270]{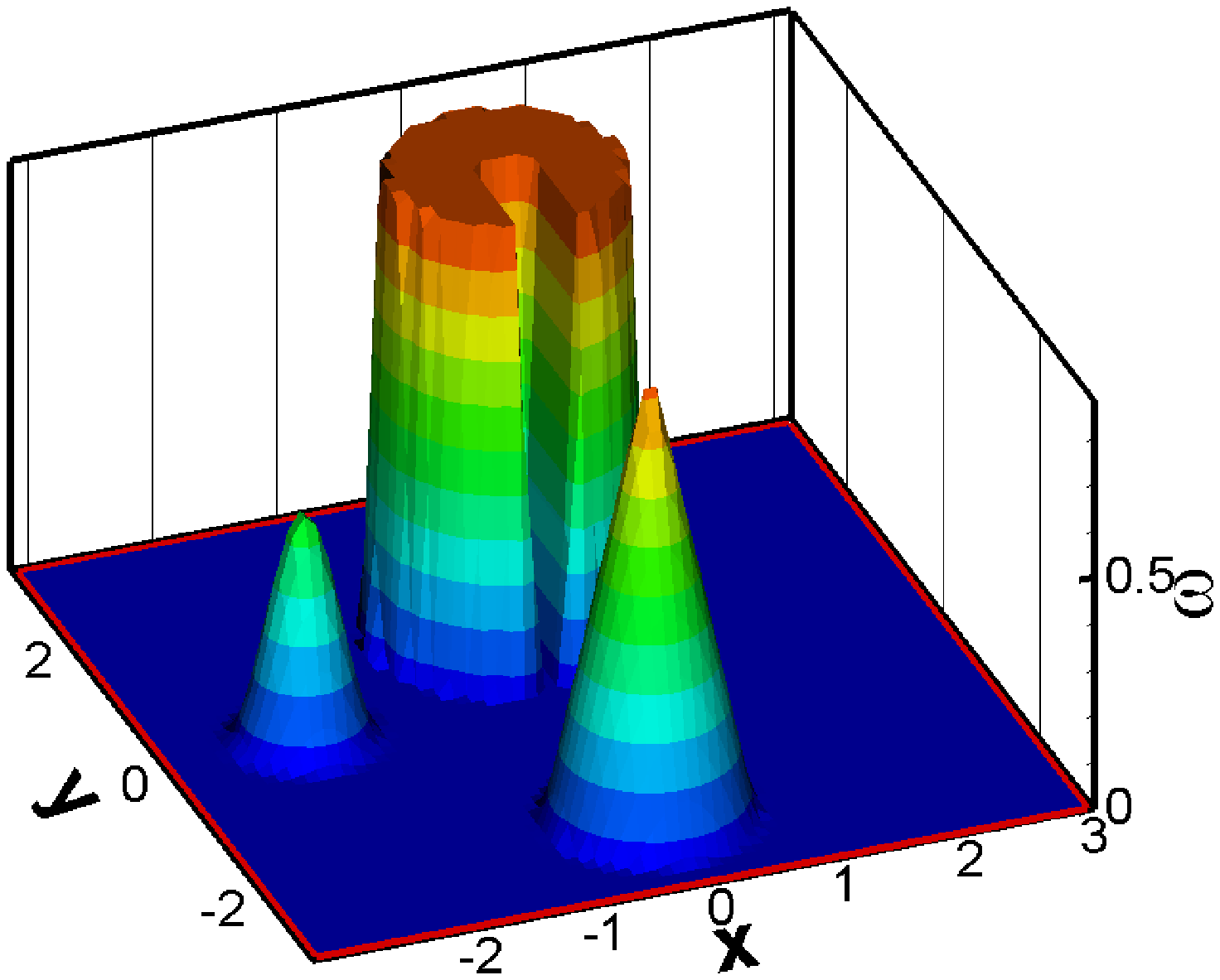}
\includegraphics[scale=0.3, angle=270]{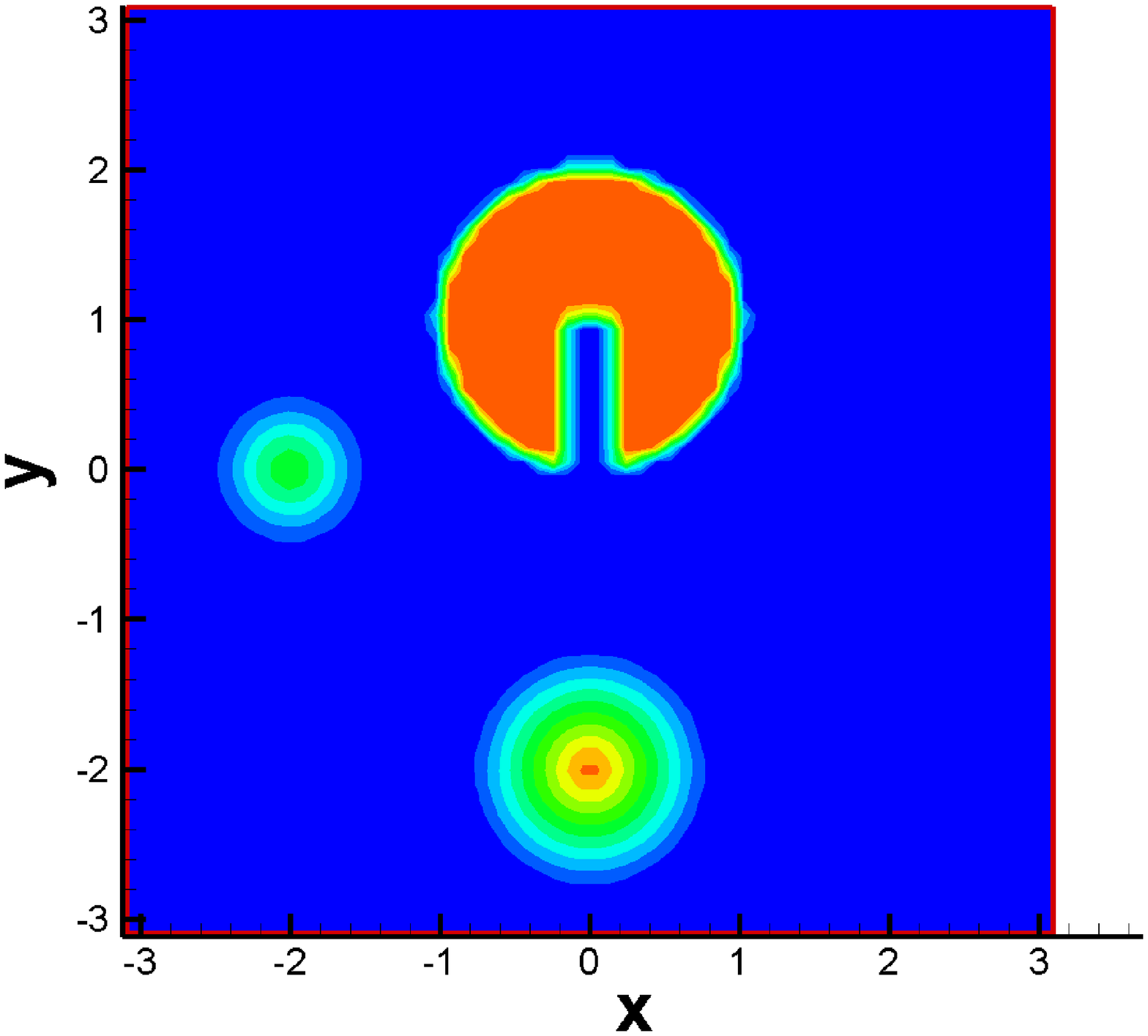} 
\end{center}
\caption{Initial profile for rigid body rotation problem (\ref{eq412}) with mesh $N^2=64^2$.}
\label{fig406}
\end{figure}


\begin{figure}[htbp]
\begin{center}
\includegraphics[scale=0.3]{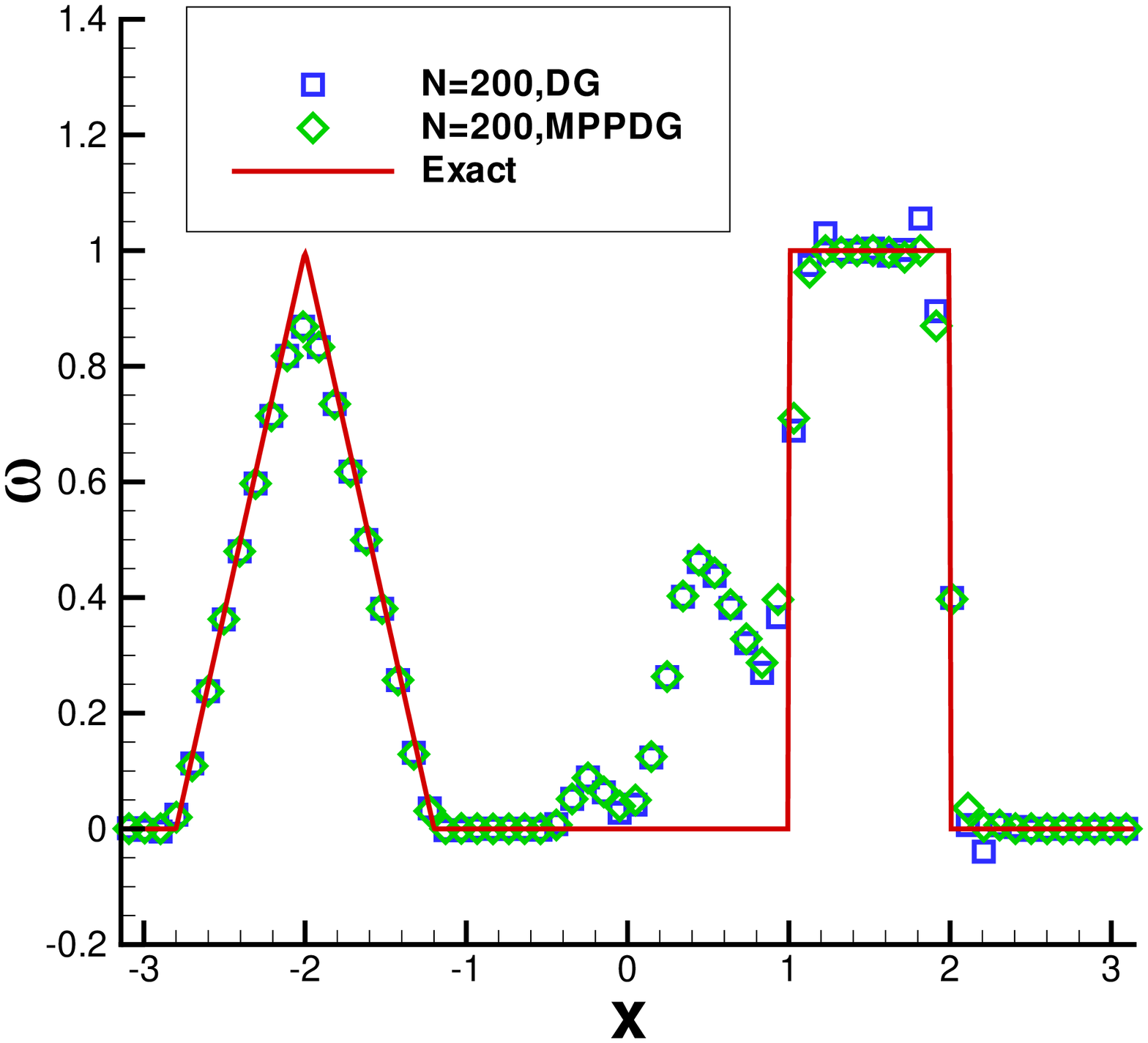},
\includegraphics[scale=0.3]{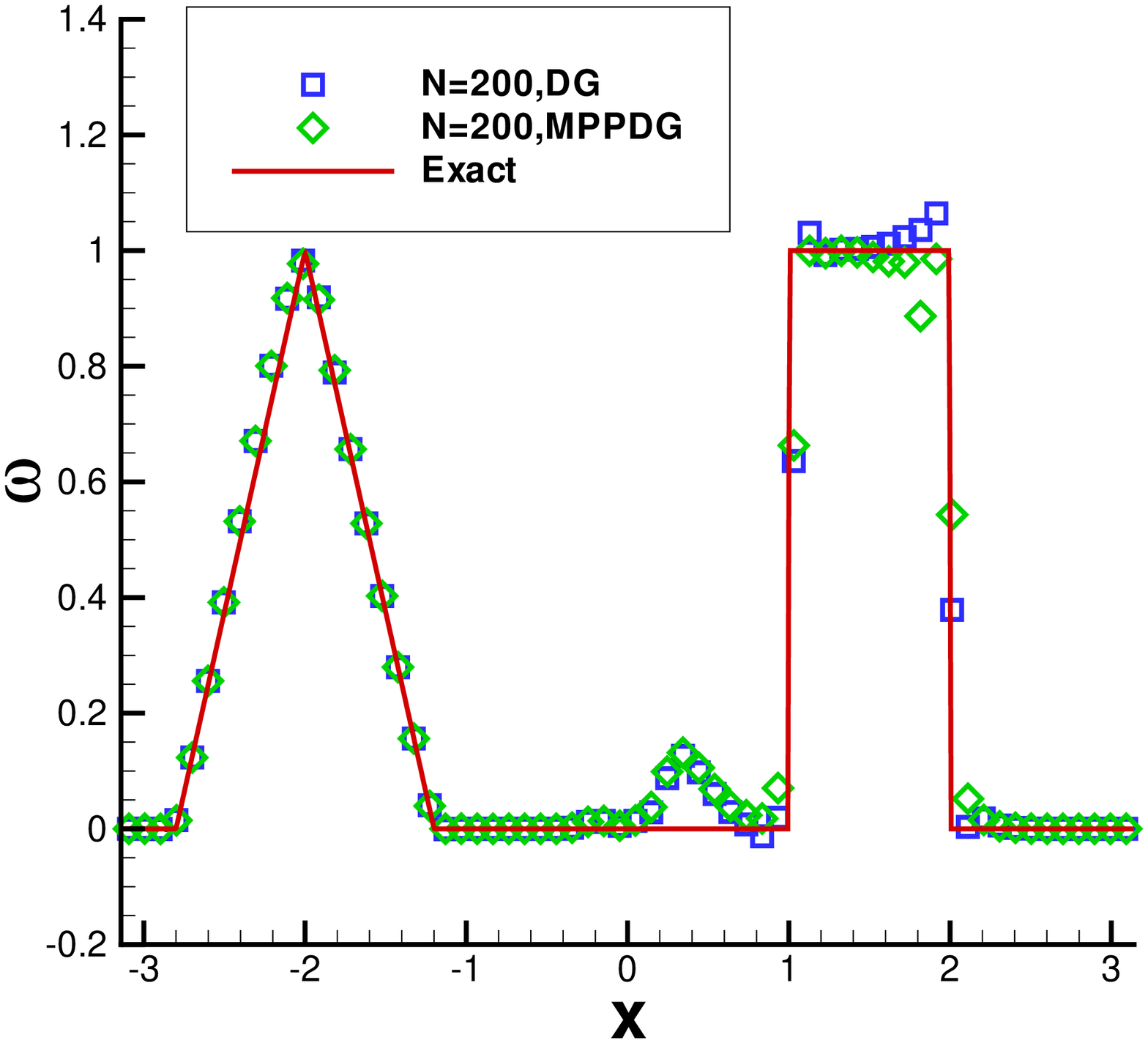} \\
\includegraphics[scale=0.3]{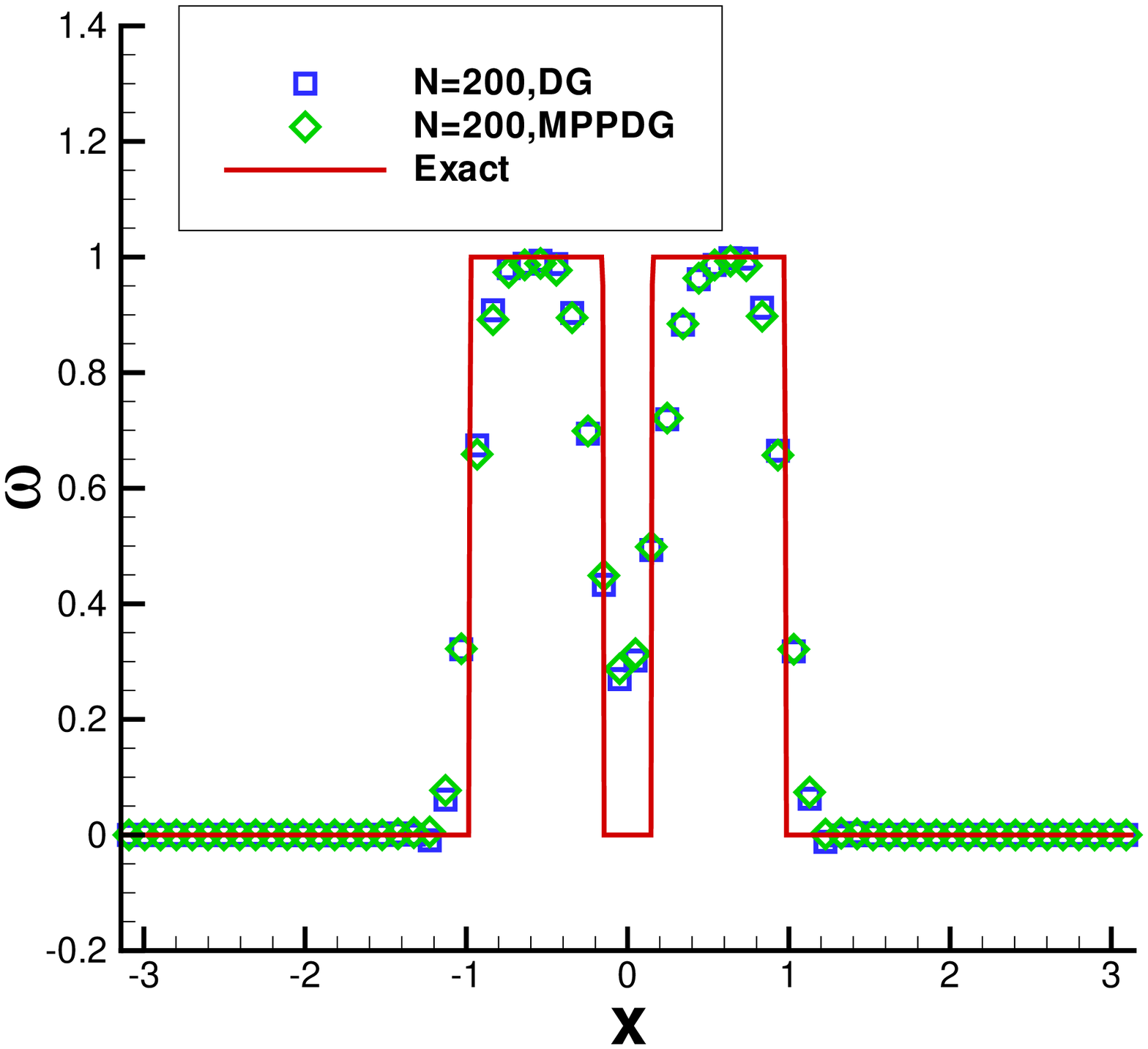},
\includegraphics[scale=0.3]{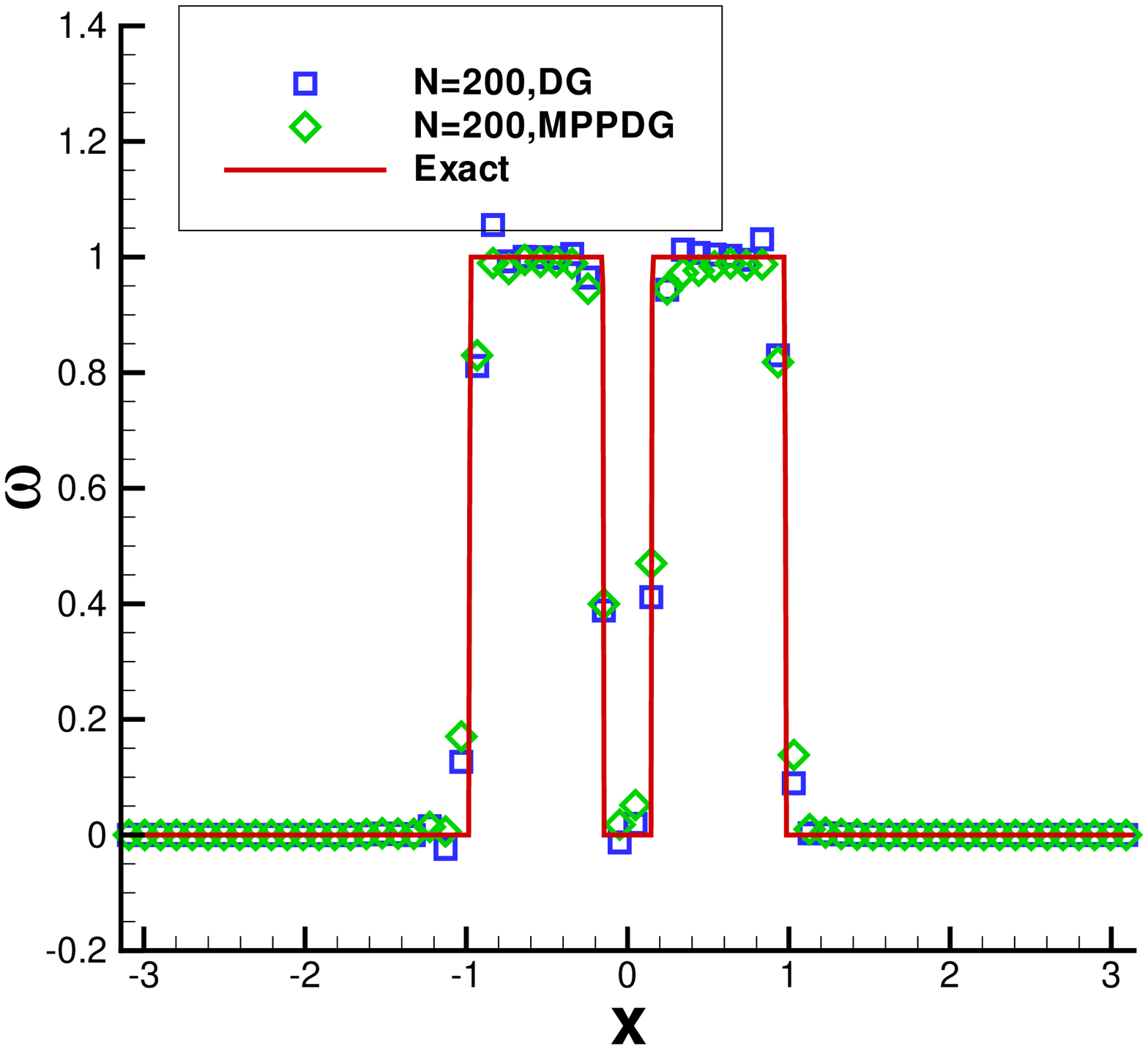} \\
\includegraphics[scale=0.3]{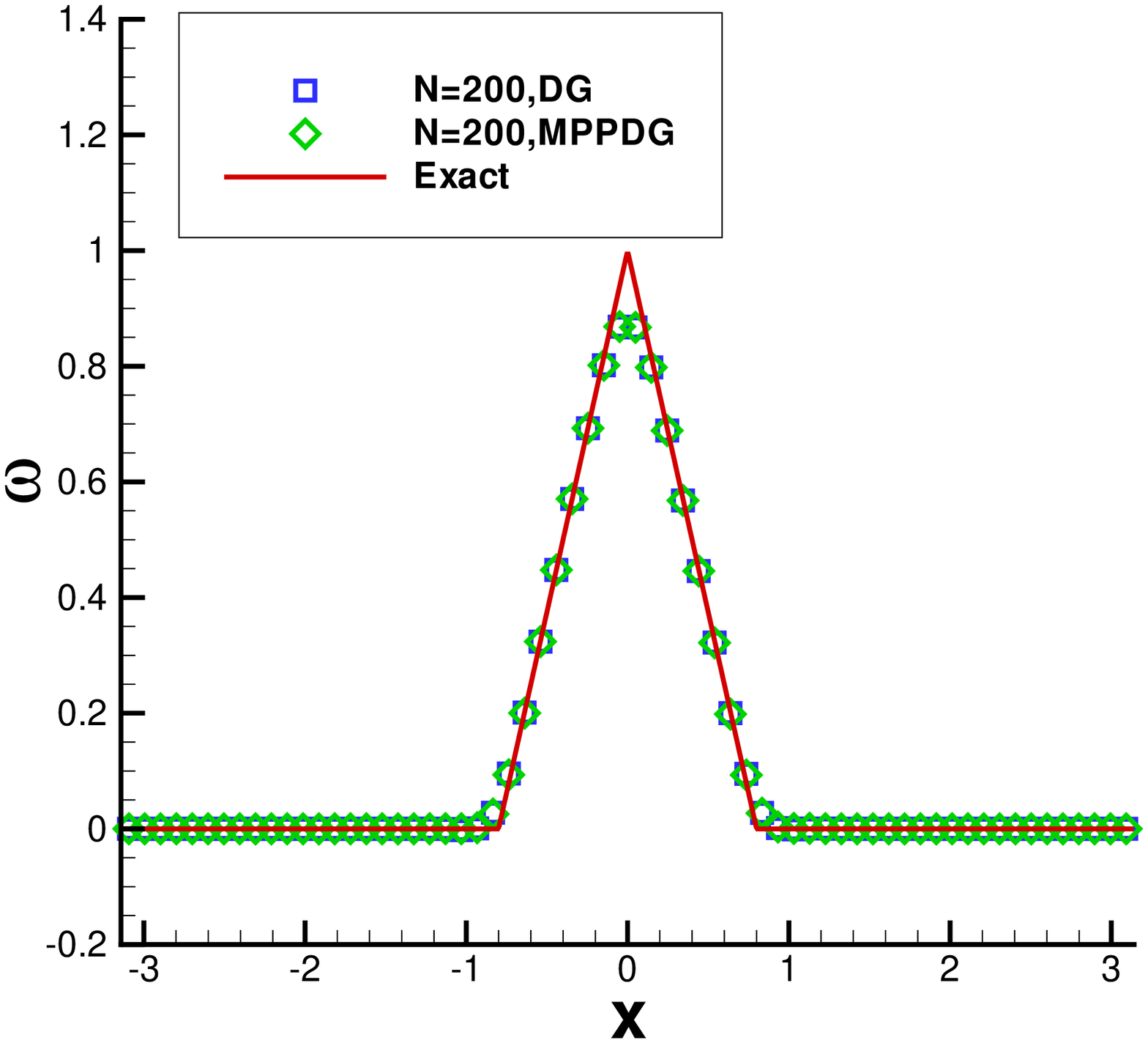},
\includegraphics[scale=0.3]{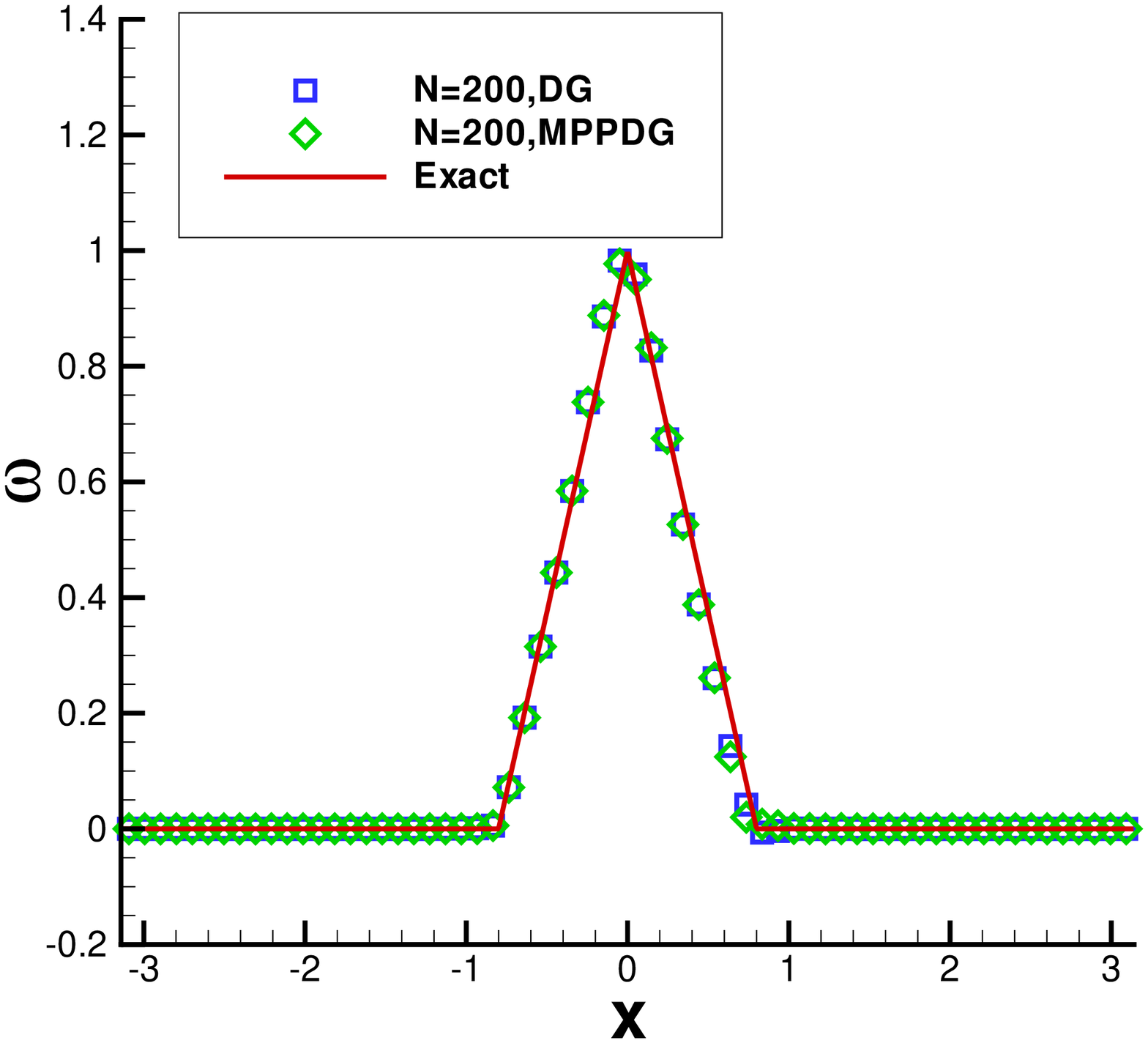}
\end{center}
\caption{Rigid body rotation problem (\ref{eq412}) (Left) and swirling deformation flow (\ref{eq413}) (Right). $T=2\pi$. $P^2$ with mesh $N^2=64^2$. From top to bottom, cuts along $x=0$, $y=0.8$ and $y=-2$ respectively.}
\label{fig407}
\end{figure}

\end{exa}

\begin{exa}(Swirling deformation flow) 
We consider the viscous swirling deformation flow
\beq
\omega_t-(\cos^2(\frac{x}{2})^2 \sin(y) g(t) \omega)_x+(\sin(x) \cos^2(\frac{y}{2}) g(t) \omega)_y=\frac{1}{Re}(\omega_{xx}+\omega_{yy}),
\label{eq413}
\eeq
with periodic boundary conditions on the domain $[-\pi, \pi]^2$ and $g(t)=\cos(\pi t/T)/\pi$. The initial condition is the same as shown in Fig. \ref{fig406}. In Fig. \ref{fig407} (right), we show the cuts along $x=0$, $y=0.8$ and $y=-2$ for the $P^2$ solutions with and without the MPP flux limiter when $1/Re=0$ at $T=2\pi$ without the TVB limiter. In Table \ref{tab407}, the minimum and maximum values of the $P^2$ solutions of $1/Re=0.01$ at $T=0.1$ with the TVB limiter are displayed; it is also observed that the undershoots and overshoots can be effectively eliminated by the MPPDG method. 


\begin{table}[ht]
\centering
\caption{Minimum and maximum values of the $P^2$ solutions for the swirling deformation flow problem (\ref{eq413}) with initial condition in Fig. \ref{fig406}. $T=0.1$.}
\vspace{0.2cm}
  \begin{tabular}{|c||c c|c c|}
    \hline
        & \multicolumn{2}{c|}{DG}               &   \multicolumn{2}{c|}{MPPDG}        \\ \hline
  $N^2$ &  $(\bar \omega_h)_{min}$    &  $(\bar \omega_h)_{max}$    &  $(\bar \omega_h)_{min}$   &  $(\bar \omega_h)_{max}$   \\ \hline
  $8^2$ & -0.0006852927049  &  0.9050469531584  &  0.0000000000001 &  0.9050413650118 \\ \hline
 $16^2$ & -0.0130846231963  &  1.0085607205967  &  0.0000000000000 &  0.9999999999999 \\ \hline
 $32^2$ & -0.0091718709627  &  1.0119090109548  &  0.0000000000000 &  0.9999999999999 \\ \hline
 $64^2$ & -0.0046303690842  &  1.0040478981792  &  0.0000000000000 &  0.9999999999999 \\ \hline
$128^2$ & -0.0000082213608  &  1.0000807250115  &  0.0000000000000 &  0.9998834799957 \\ \hline
  \end{tabular}
\label{tab407}
\end{table}

\end{exa}

\begin{exa}(Accuracy test) 
Now we consider an example with an exact smooth solution to the incompressible flow 
problems (\ref{ins1})-(\ref{ins3}), which is defined on $[0,2\pi]^2$ with periodic boundary conditions.
The exact solution is given by 
\beq
\omega(x,y,t)=-2\sin(x)\sin(y)\exp(-2t/Re).
\label{eq414}
\eeq
In Table \ref{tab408}, very slight difference can be seen between DG and MPPDG solutions,
which indicates that the high order of accuracy would not be affected by the MPP flux limiter.

\begin{table}[ht]
\centering
\caption{$L^1$ and $L^\infty$ errors and orders for the incompressible flow problem with exact solution (\ref{eq414}). $P^2$ and $T=0.1$.}
\vspace{0.2cm}
  \begin{tabular}{|c||c|c|c|c|c|c|c|}
    \hline
 &   $N^2$ &  $L^1$ error  & order  & $L^\infty$ error & order  & $(\bar \omega_h)_{min}$    & $(\bar \omega_h)_{max}$   \\ \hline
\multirow{5}{*}{DG}
 &   $8^2$ &   6.25E-03 &       --&   7.16E-02 &       --& -1.6182253025700&  1.6182245357439 \\ \cline{2-8}
 &  $16^2$ &   7.77E-04 &     3.01&   9.41E-03 &     2.93& -1.8955328534064&  1.8955328450178 \\ \cline{2-8}
 &  $32^2$ &   9.31E-05 &     3.06&   1.11E-03 &     3.08& -1.9704911703283&  1.9704911702569 \\ \cline{2-8}
 &  $64^2$ &   1.11E-05 &     3.07&   1.32E-04 &     3.07& -1.9895998331451&  1.9895998330880 \\ \cline{2-8}
 & $128^2$ &   1.35E-06 &     3.04&   1.58E-05 &     3.07& -1.9944013498433&  1.9944013498428 \\ \hline
\multirow{5}{*}{MPPDG}
 &   $8^2$ &   6.25E-03 &       --&   7.16E-02 &       --& -1.6182253025700&  1.6182245357439 \\ \cline{2-8}
 &  $16^2$ &   7.77E-04 &     3.01&   9.41E-03 &     2.93& -1.8955328534063&  1.8955328450178 \\ \cline{2-8}
 &  $32^2$ &   9.31E-05 &     3.06&   1.11E-03 &     3.08& -1.9704911703283&  1.9704911702568 \\ \cline{2-8}
 &  $64^2$ &   1.11E-05 &     3.07&   1.32E-04 &     3.07& -1.9895998331451&  1.9895998330880 \\ \cline{2-8}
 & $128^2$ &   1.35E-06 &     3.04&   1.58E-05 &     3.07& -1.9944013498433&  1.9944013498426 \\ \hline
  \end{tabular}
\label{tab408}
\end{table}
\end{exa}

\begin{exa}(Vortex patch problem) 
We consider the incompressible Navier-Stokes problems (\ref{ins1})-(\ref{ins3}) with the following initial condition
\beq
\omega(x,y,0)=
\begin{cases}
-1, \quad & \f{\pi}{2}\le x \le \f{3\pi}{2}, \f{\pi}{4}\le y \le \f{3\pi}{4}; \\
 1, \quad & \f{\pi}{2}\le x \le \f{3\pi}{2}, \f{5\pi}{4}\le y \le \f{7\pi}{4}; \\
 0, \quad & \text{otherwise}.
\end{cases}
\label{eq415}
\eeq
on the domain $[0, 2\pi]^2$ and with periodic boundary conditions. From Table \ref{tab409}, we can see
that the undershoots and overshoots of the DG solutions at $T=0.1$ have been eliminated by the MPPDG solutions too. The contour plots of $P^2$ MPPDG at $T=5$ are presented in Fig. \ref{fig408}, and the DG
solutions are similar.

\begin{table}[ht]
\centering
\caption{Minimum and maximum values of $P^2$ solutions for the vortex patch problem with initial condition (\ref{eq415}). $T=0.1$.}
\vspace{0.2cm}
\begin{tabular}{|c|c c|c c|}
\hline
        &\multicolumn{2}{c|}{DG}               &  \multicolumn{2}{c|}{MPPDG}            \\ \hline
 $N^2$  &   $(\bar \omega_h)_{min}$   &   $(\bar \omega_h)_{max}$  &   $(\bar \omega_h)_{min}$    &   $(\bar \omega_h)_{max}$   \\ \hline
 $8^2$  & -0.98522553415724 & 0.98522282105428 & -0.985225534157236 & 0.985222821054276 \\ \hline
$16^2$  & -1.01012884109198 & 0.99926259257547 & -0.994575068611273 & 0.999267288756283 \\ \hline
$32^2$  & -0.99999958730449 & 0.99999960386729 & -0.999981109945692 & 0.999999614456395 \\ \hline
$64^2$  & -1.00003514031430 & 1.00003228822257 & -0.999981418283365 & 0.999998335394442 \\ \hline
\end{tabular}
\label{tab409}
\end{table}

\begin{figure}[htbp]
\begin{center}
\includegraphics[scale=0.3]{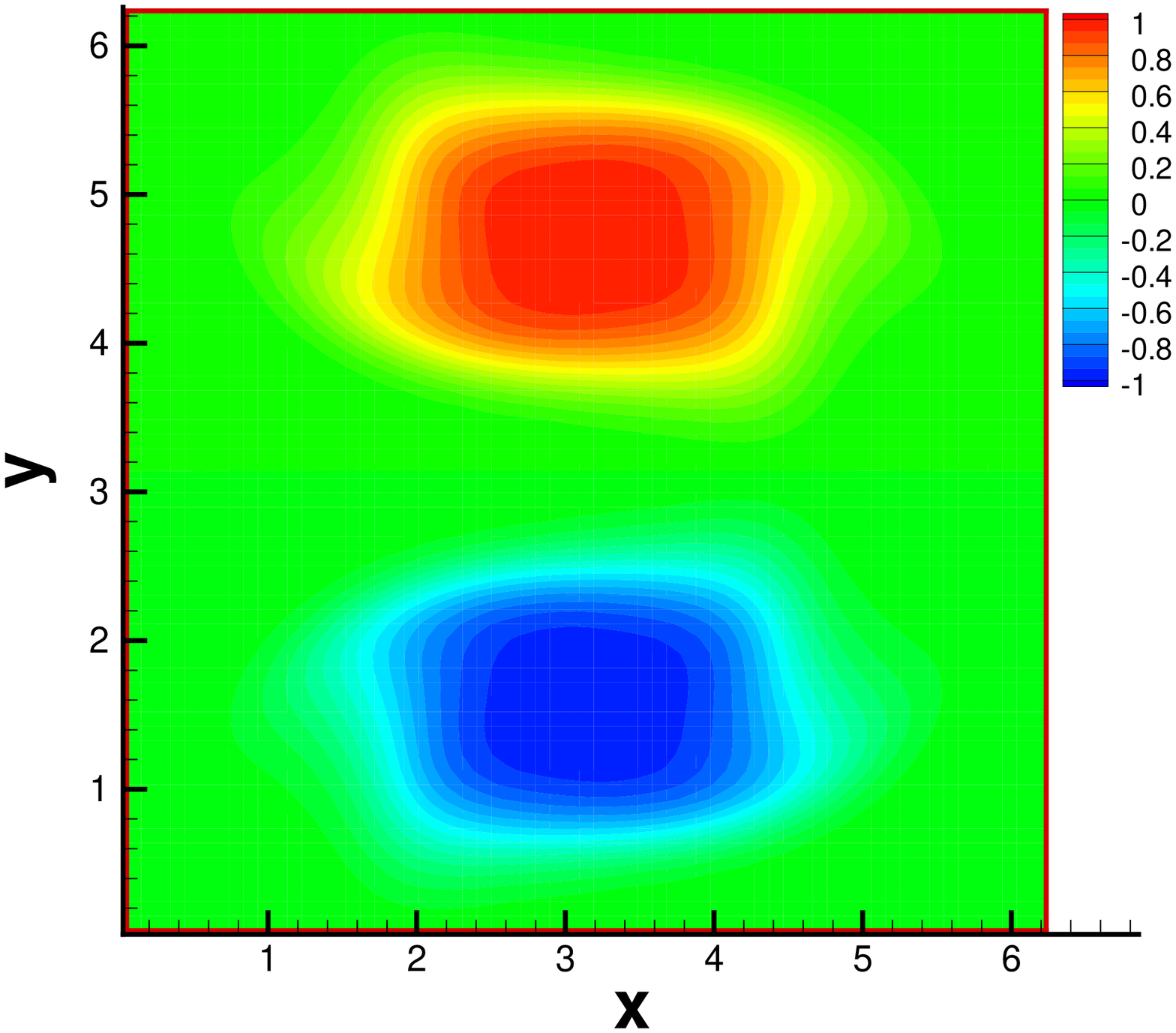},
\includegraphics[scale=0.3]{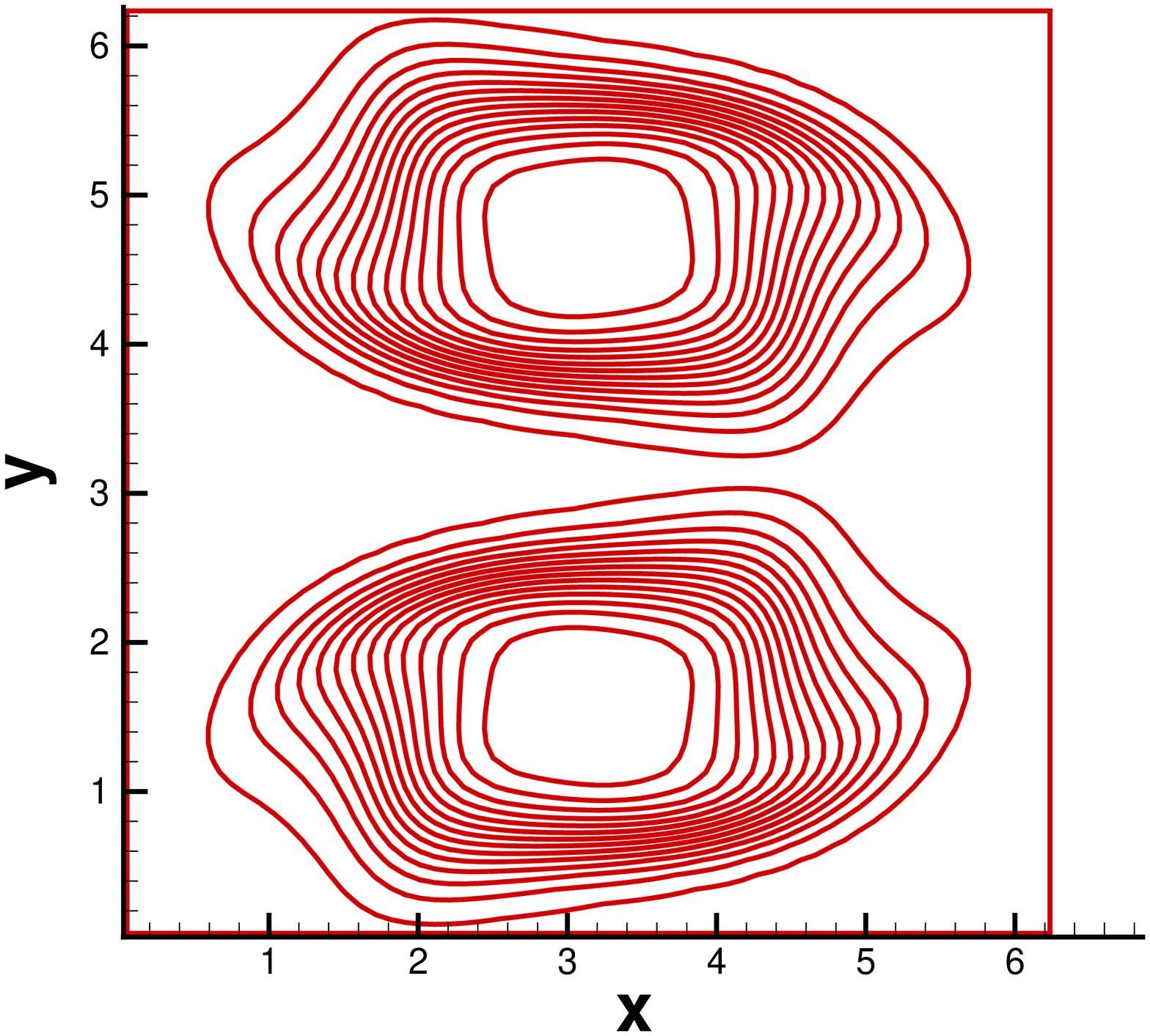}
\end{center}
\caption{Contour plots of vortex patch problem with initial condition (\ref{eq415}). $P^2$ 
MPPDG on the mesh of $N^2=64^2$. Right ones are 31 equally spaced contour lines from $-1$ to $1$. }
\label{fig408}
\end{figure} 

\end{exa}

%% file: conclusion.tex
\section{Conclusion}
\label{sec4}
\setcounter{equation}{0}
\setcounter{figure}{0}
\setcounter{table}{0}

In this paper, we propose to apply the parametrized MPP flux limiter to the RKDG method,
for solving the convection-diffusion equations. Our limiter is based on the scheme's 
conservative flux difference form in updating cell averages when the test function of the DG formulation 
is taken to be $1$. The proposed approach is mass conservative and can be applied for DG methods 
with piecewise polynomial spaces of degree $k$ ($k\ge0$).
It also has low computational cost and is easy to implement, as it is applied 
only at the final RK stage for the evolution of cell averages (not for higher moments),
in order to preserve the solution cell averages' MPP property. One major difficulty is to provide the proof of the arbitrary high order accuracy (higher than 3rd order), even though extensive numerical tests have been shown for the robust performance of the limiters in preserving 
accuracy and MPP properties of the high order numerical solutions.

%% file: log.bbl
\begin{thebibliography}{10}

\bibitem{cheng2008discontinuous}
{\sc Y.~Cheng and C.-W. Shu}, {\em {A discontinuous Galerkin finite element
  method for time dependent partial differential equations with higher order
  derivatives}}, Mathematics of Computation, 77 (2008), pp.~699--730.

\bibitem{cockburn1990runge}
{\sc B.~Cockburn, S.~Hou, and C.-W. Shu}, {\em {The {R}unge-{K}utta local
  projection discontinuous {G}alerkin finite element method for conservation
  laws. {IV}. The multidimensional case}}, Mathematics of Computation, 54
  (1990), pp.~545--581.

\bibitem{cockburn2000development}
{\sc B.~Cockburn, G.~E. Karniadakis, and C.-W. Shu}, {\em {The development of
  discontinuous Galerkin methods}}, Springer, 2000.

\bibitem{cockburn1989tvb3}
{\sc B.~Cockburn, S.-Y. Lin, and C.-W. Shu}, {\em {{TVB} {R}unge-{K}utta local
  projection discontinuous {G}alerkin finite element method for conservation
  laws {III}: one-dimensional systems}}, Journal of Computational Physics, 84
  (1989), pp.~90--113.

\bibitem{cockburn1989tvb2}
{\sc B.~Cockburn and C.-W. Shu}, {\em {{TVB} {R}unge-{K}utta local projection
  discontinuous {G}alerkin finite element method for conservation laws. {II}.
  General framework}}, Mathematics of Computation, 52 (1989), pp.~411--435.

\bibitem{cockburn1998local}
\leavevmode\vrule height 2pt depth -1.6pt width 23pt, {\em {The local
  discontinuous Galerkin method for time-dependent convection-diffusion
  systems}}, SIAM Journal on Numerical Analysis, 35 (1998), pp.~2440--2463.

\bibitem{cockburn1998runge}
\leavevmode\vrule height 2pt depth -1.6pt width 23pt, {\em {The
  {R}unge--{K}utta discontinuous {G}alerkin method for conservation laws {V}:
  multidimensional systems}}, Journal of Computational Physics, 141 (1998),
  pp.~199--224.

\bibitem{cockburn2001runge}
\leavevmode\vrule height 2pt depth -1.6pt width 23pt, {\em {Runge--Kutta
  discontinuous Galerkin methods for convection-dominated problems}}, Journal
  of Scientific Computing, 16 (2001), pp.~173--261.

\bibitem{farago2006discrete}
{\sc I.~Farag\'{o} and R.~Horv\'{a}th}, {\em {Discrete maximum principle and
  adequate discretizations of linear parabolic problems}}, SIAM Journal on
  Scientific Computing, 28 (2006), pp.~2313--2336.

\bibitem{farago2005discrete}
{\sc I.~Farag{\'o}, R.~Horv{\'a}th, and S.~Korotov}, {\em {Discrete maximum
  principle for linear parabolic problems solved on hybrid meshes}}, Applied
  Numerical Mathematics, 53 (2005), pp.~249--264.

\bibitem{farago2012discrete}
{\sc I.~Farag\'{o} and J.~Kar\'{a}tson}, {\em {Discrete maximum principle for
  nonlinear parabolic PDE systems}}, {IMA Journal of Numerical Analysis},
  (2012).

\bibitem{fujii1973some}
{\sc H.~Fujii}, {\em {Some remarks on finite element analysis of time-dependent
  field problems}}, Theory and Practice in Finite Element Structural Analysis,
  (1973), pp.~91--106.

\bibitem{Jiang_Shu}
{\sc G.-S. Jiang and C.-W. Shu}, {\em {Efficient implementation of weighted ENO
  schemes}}, Journal of Computational Physics, 126 (1996), pp.~202--228.

\bibitem{jiang2013parametrized}
{\sc Y.~Jiang and Z.~Xu}, {\em {Parametrized maximum principle preserving
  limiter for finite difference WENO schemes solving convection-dominated
  diffusion equations}}, SIAM Journal on Scientific Computing, 35 (2013),
  pp.~A2524--A2553.

\bibitem{kurganov2000new}
{\sc A.~Kurganov and E.~Tadmor}, {\em {New high-resolution central schemes for
  nonlinear conservation laws and convection--diffusion equations}}, Journal of
  Computational Physics, 160 (2000), pp.~241--282.

\bibitem{leveque2002finite}
{\sc R.~LeVeque}, {\em {Finite volume methods for hyperbolic problems}},
  vol.~31, Cambridge university press, 2002.

\bibitem{mpp_xuMD}
{\sc C.~Liang and Z.~Xu}, {\em Parametrized maximum principle preserving flux
  limiters for high order schemes solving multi-dimensional scalar hyperbolic
  conservation laws}, Journal of Scientific Computing, 58 (2014), pp.~41--60.

\bibitem{mizukami1985petrov}
{\sc A.~Mizukami and T.~J. Hughes}, {\em {A Petrov-Galerkin finite element
  method for convection-dominated flows: an accurate upwinding technique for
  satisfying the maximum principle}}, Computer Methods in Applied Mechanics and
  Engineering, 50 (1985), pp.~181--193.

\bibitem{muskat1946flow}
{\sc M.~Muskat and R.~D. Wyckoff}, {\em {The flow of homogeneous fluids through
  porous media}}, JW Edwards Ann Arbor, 1946.

\bibitem{shu1998essentially}
{\sc C.-W. Shu}, {\em {Essentially non-oscillatory and weighted essentially
  non-oscillatory schemes for hyperbolic conservation laws}}, Advanced
  numerical approximation of nonlinear hyperbolic equations,  (1998),
  pp.~325--432.

\bibitem{shu1999high}
{\sc C.-W. Shu}, {\em {High order ENO and WENO schemes for computational fluid
  dynamics}}, in High-order methods for computational physics, T.J. Barth and
  H. Deconinck, eds., Lecture Notes in Comput. Sci. Engrg. 9, Springer-Verlag,
  Berlin, 1999, pp.~439--582.

\bibitem{shu2009high}
\leavevmode\vrule height 2pt depth -1.6pt width 23pt, {\em {High order weighted
  essentially nonoscillatory schemes for convection dominated problems}}, SIAM
  review, 51 (2009), pp.~82--126.

\bibitem{shu1988efficient}
{\sc C.-W. Shu and S.~Osher}, {\em {Efficient implementation of essentially
  non-oscillatory shock-capturing schemes}}, Journal of Computational Physics,
  77 (1988), pp.~439--471.

\bibitem{vejchodsky2010discrete}
{\sc T.~Vejchodsk{\`y}, S.~Korotov, and A.~Hannukainen}, {\em {Discrete maximum
  principle for parabolic problems solved by prismatic finite elements}},
  Mathematics and Computers in Simulation, 80 (2010), pp.~1758--1770.

\bibitem{wang2013error}
{\sc H.~Wang and Q.~Zhang}, {\em Error estimate on a fully discrete local
  discontinuous galerkin method for linear convection-diffusion problem},
  Journal of Computational Mathematics, 31 (2013), pp.~283--307.

\bibitem{mpp_xqx}
{\sc T.~Xiong, J.-M. Qiu, and Z.~Xu}, {\em {A parametrized maximum principle
  preserving flux limiter for finite difference {RK-WENO} schemes with
  applications in incompressible flows }}, Journal of Computational Physics,
  252 (2013), pp.~310--331.

\bibitem{pp_euler}
\leavevmode\vrule height 2pt depth -1.6pt width 23pt, {\em {Parametrized
  positivity preserving flux limiters for the high order finite difference WENO
  scheme solving compressible Euler equations}},  (submitted).

\bibitem{mpp_xu}
{\sc Z.~Xu}, {\em {Parametrized maximum principle preserving flux limiters for
  high order scheme solving hyperbolic conservation laws: one-dimensional
  scalar problem}}, Mathematics of Computation,  (in press).

\bibitem{yang2014parametrized}
{\sc P.~Yang, T.~Xiong, J.-M. Qiu, and Z.~Xu}, {\em {High order maximum
  principle preserving finite volume method for
  convection dominated problems}},  (submitted).

\bibitem{zhang2012maximumcd}
{\sc X.~Zhang, Y.~Liu, and C.-W. Shu}, {\em Maximum-principle-satisfying high
  order finite volume weighted essentially nonoscillatory schemes for
  convection-diffusion equations}, SIAM Journal on Scientific Computing, 34
  (2012), pp.~A627--A658.

\bibitem{zhang2010maximum}
{\sc X.~Zhang and C.-W. Shu}, {\em {On maximum-principle-satisfying high order
  schemes for scalar conservation laws}}, Journal of Computational Physics, 229
  (2010), pp.~3091--3120.

\bibitem{zhang2011maximum}
\leavevmode\vrule height 2pt depth -1.6pt width 23pt, {\em
  {Maximum-principle-satisfying and positivity-preserving high-order schemes
  for conservation laws: survey and new developments}}, Proceedings of the
  Royal Society A: Mathematical, Physical and Engineering Science, 467 (2011),
  pp.~2752--2776.

\bibitem{yzhang2012maximum}
{\sc Y.~Zhang, X.~Zhang, and C.-W. Shu}, {\em {Maximum-principle-satisfying
  second order discontinuous Galerkin schemes for convection-diffusion
  equations on triangular meshes}}, Journal of Computational Physics, 234
  (2012), pp.~295--316.

\end{thebibliography}
